\documentclass[hidelinks,12pt]{amsart}

\usepackage{bbm}
\usepackage{libertine}
\usepackage{enumitem}
\usepackage[T1]{fontenc}
\usepackage{stmaryrd}
\usepackage{color}
\DeclareFontFamily{OT1}{pzc}{}
\DeclareFontShape{OT1}{pzc}{m}{it}{<-> s * [1.100] pzcmi7t}{}
\DeclareMathAlphabet{\mathpzc}{OT1}{pzc}{m}{it}
\DeclareMathAlphabet{\mathcal}{OMS}{cmsy}{m}{n}
\usepackage{fullpage}
\usepackage{mdframed}
\usepackage{graphicx}
\usepackage{calligra}
\usepackage{wasysym}
\usepackage{csquotes}
\usepackage{bm}
\usepackage{amsmath}
\usepackage[all]{xy}
\usepackage{amssymb}
\usepackage{amsthm}
\usepackage{mathrsfs}
\usepackage[OT2,T1]{fontenc}
\usepackage{centernot}
\usepackage{hyperref}

\newcommand{\DMO}[2]{\DeclareMathOperator{#1}{#2}}
\DMO{\BK}{BK}
\DMO{\FL}{FL}
\DMO{\Ann}{Ann}
\DMO{\std}{std}
\DMO{\antidiag}{antidiag}
\DMO{\locadm}{loc.adm}
\DMO{\Inj}{Inj}
\DMO{\LL}{LL}
\DMO{\Dmod}{\emph{D }-mod}
\DMO{\univ}{univ}
\DMO{\Fitt}{Fitt}
\DMO{\WD}{WD}
\DMO{\geom}{geom}
\DMO{\Fl}{Fl}
\DMO{\grad}{grad}
\DMO{\labmda}{\lambda}
\DMO{\Iw}{Iw}
\DMO{\tor}{tor}
\DMO{\coh}{coh}
\DMO{\vol}{vol}
\DMO{\semsim}{ss}
\DMO{\free}{free}
\DMO{\Alg}{Alg}
\DMO{\oth}{otherwise}
\DMO{\Ber}{Ber}
\DMO{\Diff}{Diff}
\DMO{\br}{br}
\DMO{\Isot}{Isot}
\DMO{\prim}{prim}
\DMO{\RAH}{RAH}
\DMO{\Sets}{Sets}
\DMO{\cone}{cone}
\DMO{\Grps}{Grps}
\DMO{\Dec}{Dec}
\DMO{\Flat}{Flat}
\DMO{\AbGps}{AbGps}
\DMO{\Sch}{Sch}
\DMO{\AH}{AH}
\DMO{\cl}{cl}
\DMO{\sk}{sk}
\DMO{\HC}{HC}
\DMO{\cosk}{sk}
\DMO{\ur}{ur}
\DMO{\LocSys}{LocSys}
\DMO{\rk}{rk}
\DMO{\NT}{NT}
\DMO{\cork}{cork}
\DMO{\KS}{KS}
\DMO{\MU}{MU}
\DMO{\der}{der}
\DMO{\Art}{Art}
\DMO{\Proj}{Proj}
\DMO{\End}{End}
\DMO{\Betti}{Betti}
\DMO{\Sym}{Sym}
\DMO{\cInd}{cInd}
\DMO{\GL}{GL}
\DMO{\Gal}{Gal}
\DMO{\Br}{Br}
\DMO{\Der}{Der}
\DMO{\Sp}{Sp}
\DMO{\Tan}{Tan}
\DMO{\Spin}{Spin}
\DMO{\Var}{Var}
\DMO{\Nrd}{Nrd}
\DMO{\cusp}{cusp}
\DMO{\Mat}{Mat}
\DMO{\Isom}{Isom}
\DMO{\Stab}{Stab}
\DMO{\SO}{SO}
\DMO{\Res}{Res}
\DMO{\Lie}{Lie}
\DMO{\SU}{SU}
\DMO{\Ad}{Ad}
\DMO{\ad}{ad}
\DMO{\im}{im}
\DMO{\Frob}{Frob}
\DMO{\Fr}{Fr}
\DMO{\red}{red}
\DMO{\an}{an}
\DMO{\Pic}{Pic}
\DMO{\Tor}{Tor}
\DMO{\Hdg}{Hdg}
\DMO{\id}{id}
\DMO{\pr}{pr}
\DMO{\Mor}{Mor}
\DMO{\Ext}{Ext}
\DMO{\ML}{ML}
\DMO{\PGL}{PGL}
\DMO{\SL}{SL}
\DMO{\GU}{GU}
\DMO{\GSp}{GSp}
\DMO{\GSL}{GSL}
\DMO{\Aff}{Aff}
\DMO{\NS}{NS}
\DMO{\gr}{gr}
\DMO{\Ch}{Ch}
\DMO{\QCoh}{QCoh}
\DMO{\Coh}{Coh}
\DMO{\inv}{inv}
\DMO{\Gr}{Gr}
\DMO{\Bun}{Bun}
\DMO{\Hk}{Hk}
\DMO{\GH}{GH}
\DMO{\HT}{HT}
\DMO{\LT}{LT}
\DMO{\Int}{Int}
\DMO{\UU}{U}
\DMO{\OO}{O}
\DMO{\Loc}{Loc}
\DMO{\Conn}{Conn}
\DMO{\sing}{sing}
\DMO{\si}{si}
\DMO{\Sen}{Sen}
\DMO{\MaxSpec}{MaxSpec}
\DMO{\ran}{ran}
\DMO{\coker}{coker}
\DMO{\DIV}{div}
\DMO{\Cl}{Cl}
\DMO{\Frac}{Frac}
\DMO{\VEC}{Vec}
\DMO{\Weil}{Weil}
\DMO{\SPLIT}{split}
\DMO{\Tr}{Tr}
\DMO{\val}{val}
\DMO{\pv}{p.v.}
\DMO{\disc}{disc}
\DMO{\trdeg}{tr.deg}
\DMO{\rad}{rad}
\DMO{\codim}{codim}
\DMO{\dist}{dist}
\DMO{\length}{length}
\DMO{\diam}{diam}
\DMO{\Supp}{Supp}
\DMO{\Ass}{Ass}
\DMO{\ord}{ord}
\DMO{\RE}{Re}
\DMO{\Sh}{Sh}
\DMO{\IM}{Im}
\DMO{\Tot}{Tot}
\DMO{\Bl}{Bl}
\DMO{\lcm}{lcm}
\DMO{\ann}{ann}
\DMO{\arcsinh}{arcsinh}
\DMO{\CHAR}{char}
\DMO{\MOD}{mod}
\DMO{\BB}{BB}
\DMO{\new}{new}
\DMO{\alg}{alg}
\DMO{\Irr}{Irr}
\DMO{\res}{res}
\DMO{\rank}{rank}
\DMO{\naive}{naive}
\DMO{\tors}{tors}
\DMO{\Perf}{Perf}
\DMO{\Sht}{Sht}
\DMO{\Perv}{Perv}
\DMO{\soc}{soc}
\DMO{\Mod}{Mod}
\DMO{\cyc}{cyc}
\DMO{\SC}{sc}
\DMO{\SP}{sp}
\DMO{\Deck}{Deck}
\DMO{\PSL}{PSL}
\DMO{\Area}{Area}
\DMO{\Cont}{Cont}
\DMO{\sgn}{sgn}
\DMO{\Cat}{Cat}
\DMO{\Cov}{Cov}
\DMO{\rig}{rig}
\DMO{\FSch}{FSch}
\DMO{\Rig}{Rig}
\DMO{\Spv}{Spv}
\DMO{\Spa}{Spa}
\DMO{\trace}{trace}
\DMO{\cont}{cont}
\DMO{\aff}{aff}
\DMO{\cor}{cor}
\DMO{\CH}{CH}
\DMO{\Spec}{Spec}
\DMO{\rec}{rec}
\DMO{\LGC}{LGC}
\DMO{\un}{un}
\DMO{\conj}{conj}
\DMO{\Eval}{Eval}
\DMO{\JH}{JH}
\DMO{\can}{can}
\DMO{\Fss}{Fss}
\DMO{\Speh}{Speh}
\DMO{\Ind}{Ind}
\DMO{\ch}{ch}
\DMO{\nr}{nr}
\DMO{\Swan}{Swan}
\DMO{\St}{St}
\DMO{\Ho}{Ho}
\DMO{\HH}{HH}
\DMO{\trop}{trop}
\DMO{\Jac}{Jac}
\DMO{\vir}{vir}
\DMO{\coll}{coll}
\DMO{\reg}{reg}
\DMO{\dlog}{dlog}
\DMO{\Div}{Div}
\DMO{\ab}{ab}
\DMO{\Tam}{Tam}
\DMO{\Ran}{Ran}
\DMO{\IC}{IC}
\DMO{\Sat}{Sat}
\DMO{\Rat}{Rat}
\DMO{\loc}{loc}
\DMO{\ev}{ev}
\DMO{\st}{st}
\DMO{\pst}{pst}
\DMO{\Fil}{Fil}
\DMO{\cris}{cris}
\DMO{\dR}{dR}
\DMO{\Rep}{Rep}
\DMO{\Sel}{Sel}
\DMO{\spec}{spec}
\DMO{\Spf}{Spf}
\DMO{\JL}{JL}
\DMO{\BGL}{BGL}
\DMO{\Arc}{Arc}
\DMO{\MHS}{MHS}
\DMO{\Nm}{Nm}
\DMO{\holim}{holim}
\DMO{\nInd}{nInd}
\DMO{\sSets}{s\textbf{Sets}}
\DMO{\sArt}{s\textbf{Art}}
\DMO{\BDJ}{BDJ}
\DMO{\GV}{GV}
\DMO{\BM}{BM}
\DMO{\Ord}{Ord}
\DMO{\mult}{mult}
\DMO{\WDRep}{WDRep}
\DMO{\Aut}{Aut}
\DMO{\Hom}{Hom}
\DMO{\sph}{sph}
\DMO{\Def}{Def}
\DMO{\GO}{GO}
\DMO{\diag}{diag}
\DMO{\cond}{cond}
\DMO{\ind}{ind}
\DMO{\irr}{irr}
\DMO{\RHom}{RHom}
\DMO{\sm}{sm}
\DMO{\sss}{ss}
\DMO{\sHom}{sHom}
\DMO{\Tran}{Tran}
\DMO{\Rees}{Rees}
\DMO{\lcv}{lcv} 
\DMO{\SN}{SN}
\DMO{\triv}{triv}
\DMO{\height}{ht}
\DMO{\proj}{proj}
\DMO{\Fun}{Fun}
\DMO{\cts}{cts}
\DMO{\Obj}{Obj}
\DMO{\Sing}{Sing}
\DMO{\Pro}{Pro}
\DMO{\Ig}{Ig}
\DMO{\Ha}{Ha}
\DMO{\BC}{BC}
\DMO{\RZ}{RZ}
\DMO{\supp}{supp}
\DMO{\projdim}{proj.dim}
\DMO{\Zar}{Zar}
\DMO{\Ban}{Ban}
\DMO{\LA}{LA}
\DMO{\ess}{ess}
\DMO{\op}{op}
\DMO{\Func}{Func}
\DMO{\Born}{Born}
\DMO{\Comm}{Comm}
\DMO{\Dr}{Dr}
\DMO{\LC}{LC}
\DMO{\nind}{n-ind}
\DMO{\perf}{perf}
\DMO{\charpoly}{char.poly}

\makeatletter
\def\thmhead@plain#1#2#3{%
  \thmname{#1}\thmnumber{\@ifnotempty{#1}{ }\@upn{#2}}%
  \thmnote{ {\the\thm@notefont#3}}}
\let\thmhead\thmhead@plain
\makeatother

\newtheorem*{thm1*}{Theorem}
\newtheorem*{lemma*}{Lemma}
\newtheorem*{defn1*}{Definition}

\newtheorem{thm1}{Theorem}

\newtheorem{thm2}{Theorem}[section]
\newtheorem{lem2}[thm2]{Lemma}

\newtheorem{prop1}{Proposition}

\newtheorem{coro2}[thm2]{Corollary}

\newtheorem*{prop*}{Proposition}

\newtheorem{prop2}[thm2]{Proposition}

\newtheorem*{conj*}{Conjecture}

\theoremstyle{definition}

\newtheorem{defn2}[thm2]{Definition}

\newtheorem*{defn2*}{Definition}

\newtheorem{homework}{}
\newtheorem*{prb*}{Problem}

\newtheorem*{claim*}{Claim}

\newtheorem{rmk2}[thm2]{Remark}

\newtheorem*{rmk*}{Remark}

\newtheorem{exam2}[thm2]{Example}

\newtheoremstyle{theoremdd}
  {6pt}
  {6pt}
  {}
  {0pt}
  {\bfseries}
  {.}
  { }
  {\thmname{#1}\thmnumber{ #2}\textnormal{\thmnote{ #3}}}
  \theoremstyle{theoremdd}

\newtheoremstyle{theoremee}
  {6pt}
  {6pt}
  {\itshape}
  {0pt}
  {\bfseries}
  {.}
  { }
  {\thmname{#1}\thmnumber{ #2}\textnormal{\thmnote{ #3}}}
  \theoremstyle{theoremee}

\newcommand{\xrar}[1]{\xrightarrow{#1}}

\newcommand{\xlar}[1]{\xleftarrow{#1}}

\newcommand{\riso}{\xrar{\sim}}
\newcommand{\liso}{\xlar{\sim}}

 \newenvironment{psmat}
  {\left(\begin{smallmatrix}}
  {\end{smallmatrix}\right)}

 \newenvironment{psmatrix}
  {\left(\begin{smallmatrix}}
  {\end{smallmatrix}\right)}
  
 \newenvironment{pmat}
  {\begin{pmatrix}}
  {\end{pmatrix}}

\newcommand{\textbox}[2] {\left\lbrace\parbox{#1em}{\center{#2}}\right\rbrace}

\newcommand{\wt}{\widetilde}
\newcommand{\wh}{\widehat}

\newcommand{\ov}{\overline}

\newcommand{\tr}{\operatorname{tr}}
\newcommand{\rar}{\rightarrow}

\newcommand{\lrar}{\leftrightarrow}

\newcommand{\ncom}[1]{\newcommand{#1}}
\ncom{\sbuset}{\subset}

\newcommand{\hrar}{\hookrightarrow}
\newcommand{\thrar}{\twoheadrightarrow}

\newcommand{\emphC}[1]{\textsf{\textbf{#1}}}

\DeclareSymbolFont{cyrletters}{OT2}{wncyr}{m}{n}
\DeclareMathSymbol{\Sha}{\mathalpha}{cyrletters}{"58}
\DeclareMathSymbol{\CyrE}{\mathalpha}{cyrletters}{"03}

\newcommand{\RelSpec}{\mathpzc{Spec}}	
		
\newcommand{\Del}{\nabla}

\newcommand{\bs}{\backslash}
\newcommand{\et}{\operatorname{\acute{e}t}}

\newcommand{\proet}{\operatorname{pro\acute{e}t}}

\DMO{\bmr}{\mathbbm{r}}
\DMO{\bmf}{\mathbbm{f}}
\DMO{\bmx}{\mathbbm{x}}

\newcommand{\bfone}{\mathbf{1}}

\newcommand{\bB}{\mathbb{B}}
\newcommand{\bC}{\mathbb{C}}
\newcommand{\bD}{\mathbb{D}}

\newcommand{\bF}{\mathbb{F}}
\newcommand{\bG}{\mathbb{G}}
\newcommand{\bH}{\mathbb{H}}

\newcommand{\bL}{\mathbb{L}}

\newcommand{\bN}{\mathbb{N}}

\newcommand{\bP}{\mathbb{P}}
\newcommand{\bQ}{\mathbb{Q}}
\newcommand{\bR}{\mathbb{R}}

\newcommand{\bZ}{\mathbb{Z}}

\newcommand{\bfD}{\mathbf{D}}

\newcommand{\sE}{\mathscr{E}}

\newcommand{\sH}{\mathscr{H}}

\newcommand{\cA}{\mathcal{A}}

\newcommand{\cE}{\mathcal{E}}

\newcommand{\cH}{\mathcal{H}}

\newcommand{\cL}{\mathcal{L}}
\newcommand{\cM}{\mathcal{M}}

\newcommand{\cO}{\mathcal{O}}
\newcommand{\cP}{\mathcal{P}}
\newcommand{\cQ}{\mathcal{Q}}

\newcommand{\cX}{\mathcal{X}}

\newcommand{\fY}{\mathfrak{Y}}

\newcommand{\fk}{\mathfrak{k}}

\newcommand{\fp}{\mathfrak{p}}

\newcommand{\MSRI}{\let\thefootnote\relax\footnotetext{A part of the material is based upon work supported by the National Science Foundation under Grant No. DMS-1928930 while the author was in residence at the Mathematical Sciences Research Institute in Berkeley, California during the Spring 2023 semester.}}

\newcommand{\tR}{\operatorname{R}}

\DMO{\GM}{GM}
\DMO{\coarse}{coarse}
\DMO{\str}{str}
\DMO{\topo}{top}
\DMO{\rel}{rel}
\DMO{\HIG}{HIG} 
\DMO{\Crys}{Crys}
\DMO{\Map}{Map}
\DMO{\Dol}{Dol}
\DMO{\Vect}{Vect}
\DMO{\PGSp}{PGSp}
\DMO{\gon}{gon}
\DMO{\Ssp}{sp}

\DMO{\Out}{Out}
\DMO{\comp}{comp}
\DMO{\Inv}{Inv}
\DMO{\MIC}{MIC}
\DMO{\Isoc}{Isoc}
\DMO{\la}{la}
\DMO{\corank}{corank}
\DMO{\Td}{Td}
\DMO{\hol}{hol} 
\DMO{\gen}{gen}
\begin{document}
\title{Theta characteristics and noncongruence modular forms}
\author{Gyujin Oh}\address{Department of Mathematics, Columbia University, 2990 Broadway, New York, NY 10027}\email{gyujinoh@math.columbia.edu}
\maketitle
\begin{abstract}
The Hodge bundle $\omega$ over a modular curve is a square-root of the canonical bundle twisted by the cuspidal divisor, or a \emph{theta characteristic}, due to the Kodaira--Spencer isomorphism. We prove that, in most cases, a section of a theta characteristic $\nu$ (or any odd power of it) different from $\omega$ is a \emph{noncongruence modular form}. On the other hand, we show how $\nu\ne\omega$ gives rise to a ``twisted'' analogue of the diagonal period map to a Siegel threefold, whose difference attributes to the stackiness of the moduli of abelian surfaces $\cA_{2}$. Some questions on the Brill--Noether theory of the modular curves are answered.
\end{abstract}

\tableofcontents
\section{Introduction}

The theory of theta characteristics was initiated in hopes to clarify the formulas satisfied by theta functions. Recall that a \emph{theta characteristic} of a connected compact Riemann surface is a divisor class $\Theta$ where $2\Theta$ is the canonical class. In general, aside from the parity of a theta characteristic, it is difficult to distinguish one theta characteristic from another. For example, by \cite{Atiyah}, choosing a theta characteristic of a compact complex manifold amounts to choosing a spin structure on the manifold. The set of all theta characteristics of a connected compact Riemann surface $C$ is a homogeneous space under the action of the group of $2$-torsions of the Jacobian of $C$, but there is no good choice of a basepoint.

However, there are certain curves over which there is a ``preferred choice'' of a theta characteristic. A prominent example is the case of a modular curve, which is the moduli space of elliptic curves with certain structures. If we denote the universal elliptic curve as $f:\cE\rar Y$, then the \emph{Hodge bundle} $f_{*}\Omega_{\cE/Y}^{1}$, often denoted as $\omega$, satisfies the Kodaira--Spencer isomorphism, which says that $\omega$ is a theta characteristic of $Y$ (see Theorem \ref{KSThm}). Starting off with this observation, this paper aims to investigate the following question.
\[\text{How ``special'' is $\omega$ compared to the other theta characteristics of the modular curve?}\]

\subsubsection*{Noncongruence modular forms and geometric local systems}
\hfill

The above vague question may be interpreted in various ways. Indeed, the modular curve $Y$ and its Hodge bundle $\omega$ is of fundamental importance in the study of arithmetic of modular forms, which may deem $Y$ and $\omega$ special. On the other hand, we give a much more precise result as to why $\omega$ is the most arithmetically meaningful theta characteristic of the modular curve.\footnote{What we actually study are the logarithmic versions of theta characteristics over a compact modular curve. For this purpose, we require logarithmic generalizations of various results (such as complex/$p$-adic nonabelian Hodge correspondence) in the paper. In the Introduction, we suppress this issue for the sake of exposition.}
\begin{thm1}\label{Theorem1}Let $\nu$ be a theta characteristic different from the Hodge bundle $\omega$. For $k\ge1$ odd, the sections of $\nu^{\otimes k}$ are \emph{noncongruence modular forms}.
\end{thm1}
This is Theorem \ref{CongruenceTheta} of the paper. Recall that a \emph{noncongruence modular form} is a holomorphic function over the upper half plane satisfying the properties of the modular forms, except that the level group is a finite index subgroup of $\SL_{2}(\bZ)$ that is not a congruence subgroup. As the sections of a power of $\omega$ are modular forms with the level being the level of the modular curve, this picks out $\omega$ as the unique theta characteristic of the modular curve whose sections are congruence modular forms.

We establish Theorem \ref{Theorem1} by studying the finite index subgroup $\Gamma_{\nu}\subset\SL_{2}(\bZ)$ associated to each theta characteristic $\nu$ via the nonabelian Hodge correspondence. It was already noticed by Simpson in \cite{Simpson} that a theta characteristic $L$ on a curve $C$ can be used to define a Higgs bundle $(L\oplus L^{-1},\theta)$ which corresponds via nonabelian Hodge correspondence to the monodromy representation underlying the variation of the Hodge structures induced from the complex uniformization of the curve $C$. As a section of $\nu^{\otimes k}$ is a weight $k$ modular form with level $\Gamma_{\nu}$, the content of Theorem \ref{CongruenceTheta} is that $\Gamma_{\nu}$ is \emph{not} a congruence subgroup for $\nu\ne\omega$, which we prove by slightly generalizing \cite{Kiming1}.

The difference between $\Gamma_{\nu}$ and the level of the modular curve $\Gamma:=\Gamma_{\omega}$ is mild, as $\Gamma_{\nu}$ and $\Gamma$ have the same image in $\PSL_{2}(\bZ)$, or equivalently, $\pm\Gamma_{\nu}=\pm\Gamma$. In fact, even though $\Gamma_{\nu}$ for $\nu\ne\omega$ is not a congruence subgroup, the standard representation $\rho_{\nu}:\Gamma_{\nu}\rar\GL_{2}(\bC)$ defines a local system of the modular curve which come from geometry. More precisely, in \S6, we construct a family of \emph{abelian surfaces} $W_{\nu}$, named the \emph{twisted Kuga--Sato variety}, over the modular curve, which gives rise to a geometric local system which contains $\rho_{\nu}$ as a sub-local system. 

The twisted Kuga--Sato variety $W_{\nu}$ over $Y$ gives rise to a \emph{twisted period map} $\pi_{\nu}:Y\rar\cA_{2}$ to the moduli stack of abelian surfaces. This differs from the usual diagonal period map $\pi_{\diag}:Y\rar\cA_{2}$ induced from the diagonal embedding $\bH\rar\bH\times\bH\rar\bH_{2}$ from the upper half plane to the Siegel upper half space of degree $2$, but in a very subtle way: $\pi_{\nu}\ne\pi_{\diag}$ but $p_{\nu}=p_{\diag}$, where $p_{\nu},p_{\diag}:Y\rar A_{2}$ are the corresponding period map into the coarse moduli scheme $A_{2}$ of abelian surfaces. Namely, the difference between the two period maps comes from the stacky nature of $\cA_{2}$. In \S7, we also construct a certain non-standard level structure on the Siegel modular threefold over which we see the difference between the two period maps on the level of schemes.

The twisted Kuga--Sato variety $W_{\nu}$ can be seen as realizing the ``geometric local system'' corresponding to the Higgs bundle $(E_{\nu}:=\nu\oplus\nu^{-1},\theta_{\nu})$. In \S6, we construct a variation of Hodge structures $\rho_{\nu,H}$ and a de Rham $\bZ_{p}$-local system $\rho_{\nu,p}$ on $Y$, constructed as a part of the relative $H^{1}$ of the twisted Kuga--Sato variety $W_{\nu}/Y$. These geometric local systems satisfy the following.
\begin{thm1}\label{Theorem2}The variation of Hodge structures $\rho_{\nu,H}$ is the unique variation of Hodge structures whose associated graded is $(E_{\nu},\theta_{\nu})$. The de Rham $\bZ_{p}$-local system $\rho_{\nu,p}$ restricted to $Y_{\bQ_{p}^{\nr}}$ (here, $\bQ_{p}^{\nr}$ is the maximal unramified extension of $\bQ_{p}$) is a crystalline $\bZ_{p}$-local system, and it is associated to a unique filtered convergent $F$-isocrystal whose associated graded is $(E_{\nu},\theta_{\nu})$.
\end{thm1}
The constructions of the local systems are in Definition \ref{LocSysConstruction}, and Theorem \ref{Theorem2} is a combination of Theorems \ref{UniquenessHodge} and \ref{UniquenessP}. The proofs use the usual nonabelian Hodge correspondence of Simpson \cite{Simpson} and the $p$-adic nonabelian Hodge correspondence established by Lan--Sheng--Zuo \cite{LSZ2} and Lan--Sheng--Yang--Zuo \cite{LSYZ}.

\subsubsection*{Brill--Noether theory of the modular curves and the Hodge bundle}
\hfill

There is another prominent avenue of research on the ``specialty'' of curves and line bundles, the \emph{Brill--Noether theory}. In the Brill--Noether theory, a line bundle $L$ over a smooth projective complex curve $C$ is considered \emph{special} if $h^{0}(L)$ is ``larger than usual.'' Furthermore, a curve $C$ is considered \emph{special} if there exists a certain line bundle with a larger than usual $h^{0}$. The notion of ``larger than usual'' is precise, as the fundamental theorems of the Brill--Noether theory (as developed by \cite{Kempf}, \cite{KleimanLaksov}, \cite{GriffithsHarris}, \cite{GiesekerPetri}) show various properties of a \emph{general curve} in $\cM_{g}$, the moduli of genus $g$ curves, and any curve or a line bundle violating these properties are deemed ``special''; for the precise definition, see Definition \ref{BNGeneral}.

Due to the special nature of the modular curves and their Hodge bundles, it is natural to guess that they may be special in the sense of Brill--Noether theory. Indeed, we show the following in \S5.
\begin{prop1}\label{Prop1}
Any modular curve of a sufficiently fine level is special in the sense of Brill--Noether theory.
\end{prop1}
This is Proposition \ref{ModCurveBN}. We think this result could be folkloric, but we were unable to find a reference to it, so we provide the proof. On the other hand, it seems that the Hodge bundle $\omega$ has no relationship with the Brill--Noether theory; we show in Examples \ref{SpecialW} and \ref{GeneralW} that sometimes the Hodge bundle $\omega$ is the theta characteristic with the most sections, while sometimes it is the theta characteristic with the fewest sections. The computation in Example \ref{GeneralW} approaches a noncongruence modular form of level $\Gamma_{\nu}$ for a theta characteristic $\nu\ne\omega$ as a square root of a weight $2$ modular form of congruence level, which is interesting in its own right.

It may still be true that the Hodge bundle $\omega$ is special in the sense of  Brill--Noether theory if the level is sufficiently fine enough. We end the Introduction with a heuristic for this.  It will be interesting to see if this heuristic can be made more precise.

Let $\Delta$ be the discriminant modular form, which is a cusp form of level $1$ and weight $12$, which vanishes exactly once at every cusp and nowhere else. If we denote the space of cusp forms of weight $k$ and level $\Gamma_{\nu}$ for a theta characteristic $\nu$ as $S_{k}(\Gamma_{\nu})$, then there is an injective map $S_{1}(\Gamma_{\nu})\xrar{\times\Delta}S_{13}(\Gamma_{\nu})$ given by multiplication by $\Delta$. The image of this map consists of the cusp forms of weight $13$ and level $\Gamma_{\nu}$ which vanish at every cusp to order $2$ or higher. Let $c_{1},\cdots,c_{m}$ be the cusps of the modular curve, and for $f\in S_{k}(\Gamma_{\nu})$, let $\vec{v}_{f}$ be the $m$-dimensional vector consisted of the first Fourier coefficients in the $q$-expansions of  $f$ at the cusps. If $f_{1},\cdots,f_{d}$ is a basis of $S_{13}(\Gamma_{\nu})$, let \[M_{\nu}:={\Bigg(}\vec{v}_{f_{1}}\quad\vec{v}_{f_{2}}\quad\cdots\quad\vec{v}_{f_{d}}{\Bigg)}.\]Then, $\dim S_{1}(\Gamma_{\nu})=\dim S_{13}(\Gamma_{\nu})-\rank M_{\nu}$. On the other hand, $\dim S_{13}(\Gamma_{\nu})$ is independent of $\nu$ by Riemann--Roch. Therefore, $\dim S_{1}(\Gamma_{\nu})$ is large if $\rank M_{\nu}$ is small, i.e. when there are more relations between $\vec{v}_{f}$'s, or when there are more relations between the Fourier expansions of the same modular form at different cusps. 

If $\nu=\omega$, then there are various relations between the Fourier expansions of a Hecke eigenform at different cusps, as there are Hecke operators; this principle was already used in the calculation of Example \ref{GeneralW}. On the other hand, if $\nu\ne\omega$, then $S_{13}(\Gamma_{\nu})$ is entirely consisted of noncongruence modular forms by Theorem \ref{Theorem1}, so the Hecke operators simply vanish by the result of Berger \cite{Berger}. Therefore, we expect that there will be less relation between the Fourier expansions of a noncongruence modular form at different cusps. This heuristic says that $\rank M_{\omega}$ has more reasons to be smaller than $\rank M_{\nu}$ for $\nu\ne\omega$, which converts to that $\dim S_{1}(\Gamma_{\omega})$ has more reasons to be larger than $\dim S_{1}(\Gamma_{\nu})$.
\subsection{Acknowledgements}
We would especially like to thank Isabel Vogt for the discussions on the Clifford index which influenced our work at its initial stage. We also thank Nathan Chen, Maarten Derickx, Filip Najman, Petar Orli\'c, Koji Shimizu, John Voight, and Tonghai Yang for helpful discussions. A part of this material is based upon work supported by the National Science Foundation under Grant No. DMS-1928930 while the author was in residence at the Mathematical Sciences Research Institute (MSRI/SLMath) in Berkeley, California, during the Spring 2023 semester.
\subsection{Notations}
Let $\Gamma\le\SL_{2}(\bZ)$ be a congruence subgroup that satisfies the following condition\footnote{We impose this condition just for simplicity, and we expect our results to be extended to more general torsion-free congruence subgroups. On the other hand, the torsion-free-ness is a more crucial assumption.}.
\begin{equation}\parbox{30em}{There exist integers $N_{1},N_{2}$ such that $(N_{1},N_{2})$ is odd, $\lcm(N_{1},N_{2})\ge5$, and $\Gamma=\Gamma_{1}(N_{1})\cap\Gamma(N_{2})$.} \tag{$*$}
\end{equation}For example, the standard congruence subgroups $\Gamma_{1}(N)$ and $\Gamma(N)$ for any $N\ge5$ satisfy ($*$). Note that ($*$) implies that $\Gamma$ is torsion-free. 

Let $Y(\Gamma)=\Gamma\bs\bH$ be the (open) modular curve, regarded as a Riemann surface, and let $X(\Gamma)$ be the compactification of $Y(\Gamma)$. Thanks to ($*$), there is a universal elliptic curve $f:\cE\rar Y(\Gamma)$. Let $D=X(\Gamma)-Y(\Gamma)$ be the cuspidal divisor. We will add subscripts to these geometric objects (e.g. $Y(\Gamma)_{\bQ}$) if we need to specify the base ring. 
Throughout the paper, we fix the embeddings $\ov{\bQ}\hrar\bC$ and $\ov{\bQ}\hrar\ov{\bQ}_{p}$, and an isomorphism $\bC\cong\ov{\bQ}_{p}$ compatible with the embeddings. We also fix a $\ov{\bQ}$-point $*\in Y(\Gamma)_{\bQ}(\ov{\bQ})$ which we use as the basepoint for $\pi_{1}$ throughout the paper. The points induced from $*$ by the embeddings $\ov{\bQ}\hrar\bC$ and $\ov{\bQ}\hrar\ov{\bQ}_{p}$ are again denoted as $*$ by abuse of notation.

We denote the genus of $X(\Gamma)$ as $g_{\Gamma}$ and the number of cusps as $n_{\Gamma}$, and we will omit the subscripts when there is no confusion. As $\Gamma$ is torsion-free, we have
\[g_{\Gamma}=1+\frac{[\SL_{2}(\bZ):\Gamma]}{24}-\frac{n_{\Gamma}}{2}.\]

The space of weight $k$ modular forms (cusp forms, respectively) of level $\Gamma$ is denoted as $M_{k}(\Gamma)$ ($S_{k}(\Gamma)$, respectively). 

Let $\cA_{g}$ be the moduli space of principally polarized abelian varieties of dimension $g$, regarded as a Deligne--Mumford stack over $\bQ$. More generally, for a level structure $\Gamma$, let $\cA_{g,\Gamma}$ be the corresponding moduli space with the $\Gamma$-level structure. Let $A_{g}$, $A_{g,\Gamma}$ be the associated coarse moduli schemes.
\section{The Hodge bundle $\omega$}
\begin{defn2}
For a field $F$ of characteristic away from the level, the \emph{Hodge bundle} $\omega$ is a line bundle over $Y(\Gamma)_{F}$  defined as 
\[\omega:=f_{*}\Omega_{\cE_{F}/Y(\Gamma)_{F}}^{1}.\]
The Hodge bundle extends canonically (in the sense of Deligne and Harris) over $X(\Gamma)_{F}$, and, by abuse of notations, we will also denote the canonical extension as $\omega$. One may, for example, define $\omega$ over $X(\Gamma)$ as the algebraization of the analytic sheaf of sections of logarithmic growth of $\omega$ at infinity over $Y(\Gamma)$. 
\end{defn2}
The following is well-known.
\begin{thm2}[(Kodaira--Spencer isomorphism)]\label{KSThm}
Over $Y(\Gamma)_{F}$, one has a canonical isomorphism
\[\KS:\omega^{\otimes2}\riso\Omega^{1}_{Y(\Gamma)_{F}/F}.\]
Over $X(\Gamma)_{F}$, one has a natural isomorphism
\[\KS:\omega^{\otimes2}\riso\Omega^{1}_{X(\Gamma)_{F}/F}(D).\]
\end{thm2}
\begin{proof}
There is a canonical morphism \[\omega\rar\omega^{-1}\otimes\Omega^{1}_{Y(\Gamma)/\bC},\]which is the Higgs field arising as the associated graded of the Gauss--Manin connection on $\sH^{1}_{\dR}(\cE_{F}/Y(\Gamma)_{F})$. Since the Gauss--Manin connection has no singularities on $Y(\Gamma)_{F}$, the Higgs field is nonvanishing everywhere, thus an isomorphism. The Kodaira--Spencer isomorphism over $X(\Gamma)_F$ follows by taking the canonical extension of both sides of the above isomorphism over $Y(\Gamma)_F$.
\end{proof}
It will be later important that  $\omega^{\otimes 2}$ and $\Omega_{X(\Gamma)}^{1}(D)$ are not just merely isomorphic to each other but also that there is a canonical isomorphism between the two.
\begin{coro2}The degree of $\omega$ is $g-1+\frac{n}{2}$.
\end{coro2}
Because of the condition on the level, we have
\[M_{k}(\Gamma)=H^{0}(X(\Gamma),\omega^{\otimes k}),\quad S_{k}(\Gamma)=H^{0}(X(\Gamma),\omega^{\otimes k}(-D)).\]A simple application of Riemann--Roch yields the following result.
\begin{prop2}\label{RR}\hfill
\begin{enumerate}\item
If $k\ge2$, we have $\dim M_{k}(\Gamma)=(k-1)(g-1)+\frac{nk}{2}$. \item If $k\ge3$, we have $\dim S_{k}(\Gamma)=(k-1)(g-1)+\frac{n(k-2)}{2}$. We also have $\dim S_{2}(\Gamma)=g$.
\item We have $\dim M_{1}(\Gamma)-\dim S_{1}(\Gamma)=\frac{n}{2}$. \item If $n>2g-2$, we have $\dim M_{1}(\Gamma)=\frac{n}{2}$ and $\dim S_{1}(\Gamma)=0$. 
\end{enumerate}
\end{prop2}
\begin{proof}
Only $\dim M_{1}(\Gamma)-\dim S_{1}(\Gamma)=\frac{n}{2}$ requires an explanation. From the short exact sequence $0\rar \omega(-D)\rar\omega\rar\omega|_{D}\rar0$, we have the long exact sequence
\[0\rar S_{1}(\Gamma)\rar M_{1}(\Gamma)\rar H^{0}(\omega|_{D})\rar H^{1}(\omega(-D))\rar H^{1}(\omega)\rar0,\] as $\omega|_{D}$ is a skyscraper sheaf. On the other hand, by Serre duality,
\[\ker(H^{1}(\omega(-D))\rar H^{1}(\omega))=\ker(H^{0}(\omega)^{*}\rar H^{0}(\omega(-D))^{*})\]\[=\left(\coker(H^{0}(\omega(-D))\rar H^{0}(\omega))\right)^{*}=\left(\frac{M_{1}(\Gamma)}{S_{1}(\Gamma)}\right)^{*}.\]Thus, we have a short exact sequence
\[0\rar\frac{M_{1}(\Gamma)}{S_{1}(\Gamma)}\rar H^{0}(\omega|_{D})\rar\left(\frac{M_{1}(\Gamma)}{S_{1}(\Gamma)}\right)^{*}\rar0.\]Therefore, $\dim M_{1}(\Gamma)-\dim S_{1}(\Gamma)$ is the half of $\dim H^{0}(\omega|_{D})=n$.
\end{proof}
\begin{rmk2}
For a cusp form of weight $1$ and level $\Gamma$ to exist, the inequality $n\le 2g-2$, or equivalently the inequality $24n\le[\SL_{2}(\bZ):\Gamma]$, must be satisfied, which is true when $\Gamma$ is sufficiently small. For example, if $\Gamma=\Gamma(N)$, the inequality is satisfied if $N\ge12$.
\end{rmk2}
\begin{rmk2}
It is expected that there is no simple formula that expresses $\dim S_{1}(\Gamma)$. It is however conjectured that $S_{1}(\Gamma)$ is mostly consisted of dihedral forms (for example, see \cite[\S1]{Duke}). 
\end{rmk2}
\section{Theta characteristics as uniformizing logarithmic Higgs bundles}
In the previous section, the computation of $\deg\omega$ and the dimension of the space of modular forms only used the Kodaira--Spencer isomorphism. Thus, the same dimension formulae will hold true for any line bundle $\nu$ such that $\nu^{\otimes2}\cong\Omega^{1}_{X(\Gamma)/\bC}(D)$.
\begin{defn2}
A line bundle $\nu$ over $X(\Gamma)$ which satisfies \[\nu^{\otimes2}\cong\Omega^{1}_{X(\Gamma)/\bC}(D),\]is called a \emph{theta characteristic}.\footnote{Perhaps a better terminology will be the \emph{stable} theta characteristic: a line bundle is called a theta characteristic if it is a square root of the canonical bundle $\Omega^{1}$; a stable theta characteristic is when a line bundle is a square root of the canonical bundle twisted by a specific divisor. As we will only care about the square-roots of $\Omega^{1}(D)$ in this article, most of the time we will just refer to such line bundles as theta characteristics.} 
\end{defn2}
\begin{lem2}
The results of Proposition \ref{RR} holds for any theta characteristic $\nu$, if we interpret $M_{k}=H^{0}(X(\Gamma),\nu^{\otimes k})$ and $S_{k}=H^{0}(Y(\Gamma),\nu^{\otimes k})$.
\end{lem2}
If $\nu$ is a theta characteristic, $\nu\otimes\omega^{-1}$ is a square-root of $\cO_{X(\Gamma)}$. Thus, there are in total $2^{2g}=\#\Jac(X(\Gamma))[2](\bC)$ many theta characteristics up to isomorphism. For a theta characteristic $\nu$, the isomorphism $\nu^{\otimes 2}\cong\Omega^{1}_{X(\Gamma)/\bC}(D)$ induces an isomorphism $\nu\cong \nu^{-1}\otimes\Omega^{1}_{X(\Gamma)/\bC}(D)$. This in turn deduces a \emph{logarithmic Higgs field} $\theta_{\nu}:E\rar E\otimes\Omega_{X(\Gamma)/\bC}^{1}(D)$ on the vector bundle $E_{\nu}:=\nu\oplus\nu^{-1}$,
\[\theta_{\nu}:\nu\oplus\nu^{-1}\rar\nu\riso\nu^{-1}\otimes\Omega^{1}_{X(\Gamma)/\bC}(D)\rar(\nu\oplus\nu^{-1})\otimes\Omega^{1}_{X(\Gamma)/\bC}(D),\]making $(E_{\nu},\theta_{\nu})$ a logarithmic Higgs bundle on $X(\Gamma)$. 

In view of the nonabelian Hodge correspondence, one may ask which local systems correspond to the Higgs fields constructed using theta characteristics. In the non-logarithmic setting, Simpson showed in \cite{Simpson} that the Higgs field formed by a theta characteristic of a hyperbolic curve is precisely a lift of the projective representation of the topological $\pi_{1}$ of the curve given by the complex uniformization.

Using a tame regular analogue of the nonabelian Hodge correspondence over a noncompact Riemann surface, we can show that the theta characterstics in our sense are also characterized by the projective lifts of $\pi_{1}(Y(\Gamma),*)=P\Gamma$. 

Before stating the theorem, we introduce some terminologies.
\begin{defn2}
Let $P\Gamma$ be the projective image of $\Gamma$. Namely, 
\[P\Gamma=\im(\Gamma\hrar\SL_{2}(\bR)\rar\PSL_{2}(\bR)).\]
A \emph{projective lift} of $P\Gamma$ is a subgroup $\Gamma'\le\SL_{2}(\bR)$ such that $P\Gamma'=P\Gamma$. A projective lift is \emph{honest} if the natural map $\Gamma'\rightarrow P\Gamma'=P\Gamma$ is injective (thus bijective).

As $Y(\Gamma)$ is topologically just a surface of genus $g$  with $n$ punctures, we can choose a set of generators $a_{1},b_{1},\cdots,a_{g},b_{g},c_{1},\cdots,c_{n}\in P\Gamma$ such that the only relation between the generators is \[[a_{1},b_{1}]\cdots[a_{g},b_{g}]c_{1}\cdots c_{n}=1.\]Let $A_{1},B_{1},\cdots,A_{g},B_{g},C_{1},\cdots,C_{n}\in\Gamma$ be the corresponding elements in $\Gamma$. A \emph{hyperbolic projective lift} is an honest projective lift of the form
\[\langle\epsilon_{11}A_{1},\epsilon_{12}B_{1},\cdots,\epsilon_{g1}A_{g},\epsilon_{g2}B_{g},C_{1},\cdots,C_{n}\rangle\subset\SL_{2}(\bR),\]where $\epsilon_{ij}\in\lbrace\pm1\rbrace$ for $1\le i\le g$, $1\le j\le 2$.
\end{defn2}
\begin{lem2}\hfill\begin{enumerate}\item
The notion of the hyperbolic projective lifts does not depend on the choice of a presentation of $P\Gamma$ as a topological fundamental group of a surface. \item Given a hyperbolic projective lift $\Gamma'$ of $P\Gamma$, let $\rho_{\Gamma'}$ be the two-dimensional real representation of $P\Gamma$ given by the composition \[\rho_{\Gamma'}:P\Gamma\liso\Gamma'\hrar\SL_{2}(\bR)\subset\GL_{2}(\bR).\]For two different hyperbolic projective lifts $\Gamma_{1}'\ne\Gamma_{2}'$, we have $\rho_{\Gamma_{1}'}\not\cong\rho_{\Gamma_{2}'}$.
\end{enumerate}
\end{lem2}
\begin{proof}\hfill
\begin{enumerate}
\item As $n\ge1$, $\Gamma=P\Gamma$ is a free group with $2g+n-1$ generators. Thus, given a presentation of $P\Gamma$ as above, choosing an honest projective lift of $P\Gamma$ is the same as choosing a sign for each of $A_{1},B_{1},\cdots,A_{g},B_{g},C_{1},\cdots,C_{n-1}$, or equivalently, choosing a homomorphism $P\Gamma\rar(\bZ/2\bZ)^{2g+n-1}$. 

Note that $X(\Gamma)$ has the fundamental group, denoted $\ov{P\Gamma}$, (with the same choice of basepoint as $Y(\Gamma)$ via the inclusion $Y(\Gamma)\hrar X(\Gamma)$) whose presentation can be given by
\[\ov{P\Gamma}\cong\langle \ov{a_{1}},\ov{b_{1}},\cdots,\ov{a_{g}},\ov{b_{g}}~|~[\ov{a_{1}},\ov{b_{1}}]\cdots[\ov{a_{g}},\ov{b_{g}}]=1\rangle,\]and the natural homomorphism $\pi_{1}(Y(\Gamma),*)\rar\pi_{1}(X(\Gamma),*)$ is given by $a_{i}\mapsto \ov{a_{i}}$, $b_{i}\mapsto \ov{b_{i}}$, $c_{i}\mapsto1$. Thus, an honest projective lift is a hyperbolic projective lift if and only if the corresponding homomorphism $P\Gamma\rar(\bZ/2\bZ)^{2g+n-1}$ factors through the morphism $P\Gamma\rar\ov{P\Gamma}$. Since the latter condition does not refer to a specific presentation at all and only uses the natural map $P\Gamma\rar\ov{P\Gamma}$, the notion of hyperbolic projective lifts is independent of the choice of a presentation of $P\Gamma$.
\item Choose a presentation of $P\Gamma$ as above. Given $\Gamma_{1}'\ne\Gamma_{2}'$, there is some $1\le i\le g$ such that either $a_{i}$ or $b_{i}$  is lifted to matrices with the opposite signs. Let $d\in\lbrace a_{i},b_{i}\rbrace$ be such element. Then, $\tr\rho_{\Gamma_{1}'}(d)=-\tr\rho_{\Gamma_{2}'}(d)$. Since $A_{i}$ and $B_{i}$ are hyperbolic matrices, $\tr A_{i}$ and $\tr B_{i}$ are both nonzero. Thus, $\tr\rho_{\Gamma_{1}'}(d)\ne\tr\rho_{\Gamma_{2}'}(d)$, which means that as abstract representations $\rho_{\Gamma_{1}'}$ and $\rho_{\Gamma_{2}'}$ are non-isomorphic.
\end{enumerate}
\end{proof}
The above lemma shows that we can refer to the hyperbolic projective lifts of $P\Gamma$ as being certain two-dimensional real representations of $P\Gamma$, or, after conjugation, two-dimensional representations of $P\Gamma$ valued in $\SU(1,1)$.

We are now able to state the main theorem of this section, which should be a standard consequence of the tame regular nonabelian Hodge correspondence.
\begin{thm2}[(Theta characteristics are hyperbolic projective lifts)]\label{MainTheta}
There is a  one-to-one correspondence between the theta characteristics and the hyperbolic projective lifts of $P\Gamma$, characterized as follows. 
\begin{itemize}
\item For a theta characteristic $\nu$, $\rho_{\Gamma_{\nu}}$ is the $2$-dimensional local system on $Y(\Gamma)$ corresponding to the logarithmic Higgs bundle $(E_{\nu},\theta_{\nu})$ via the tame nonabelian Hodge correspondence. Furthermore, there is a natural isomorphism $H^{0}(X(\Gamma),\nu^{\otimes k})\cong M_{k}(\Gamma_{\nu})$ for $k\ge1$.
\item $\Gamma_{\omega}=\Gamma$.
\end{itemize}
\end{thm2}
\begin{proof}
We would like to use the tame regular version of nonabelian Hodge correspondence over a noncompact curve as in \cite{SimpsonNoncompact}: for the definitions of the terms, see \cite[Synopsis]{SimpsonNoncompact}. 
\begin{thm2}[({Tame nonabelian Hodge correspondence over non-compact curves, {\cite[p.718]{SimpsonNoncompact}}})]Over a smooth algebraic noncompact curve, there is a natural one-to-one correspondence between stable filtered regular Higgs bundles of degree zero, and stable filtered local systems of degree zero. The correspondence preserves the rank on both sides.
\end{thm2}
On the other hand, a special case of this correspondence is proved in \cite[Theorem 4]{Simpson}: taking the graded pieces gives an equivalence of categories from the category of complex variations of Hodge structures on $Y(\Gamma)$ to the category of Hodge bundles on $Y(\Gamma)$. Here, the geometric objects on $Y(\Gamma)$ are extended to $X(\Gamma)$ as ``canonical extensions'' (namely, the filtration is given by the growth behavior at the punctures).

We equip the Higgs bundle $(E_{\nu},\theta_{\nu})$ with a left-continuous decreasing filtration
\[E_{\nu,\alpha}=E_{\nu}(-\lceil \alpha\rceil D),\quad\alpha\in\bR.\]This is by definition a filtered regular Higgs bundle of degree zero. Moreover, it is stable, as the only proper nonzero $\theta$-stable subbundle of $E_{\nu}$ is $\nu^{-1}$, whose filtered degree is negative. This is the same as the ``canonical extension'' of $(E_{\nu},\theta_{\nu})\rvert_{Y(\Gamma)}$.

By the tame nonabelian Hodge correspondence, from $(E_{\nu},\lbrace E_{\nu,\alpha}\rbrace,\theta_{\nu})$, we obtain a $2$-dimensional stable filtered local system $L_{\nu}$ of degree zero. The correspondence of the statement of the Theorem is then \[\nu\mapsto\text{the underlying local system of }L_{\nu}.\]
The inverse of the correspondence can be given as follows. Let $\Gamma'$ be a hyperbolic projective lift of $P\Gamma$. Then, the universal variation of Hodge structures on $\bH$ descend to a variation of Hodge structure on $Y(\Gamma')=Y(\Gamma)$ whose underlying local system is the same as the local system corresponding to $\Gamma'$. Since the local system has unipotent local monodromies around the punctures, the Hodge filtration extends canonically (in the sense of Deligne) to $X(\Gamma)$ as a filtration of vector bundles. Let $\ov{F^{1}}$ be the canonical extension of $F^{1}$; namely, it is the sheaf of sections of $F^{1}$ with at worst logarithmic growth at the punctures. Then the inverse correspondence is
\[\Gamma'\mapsto\ov{F^{1}}.\]This is certainly a restriction of the inverse of the tame nonabelian Hodge correspondence as above by \cite[Theorem 4]{Simpson}. It sends hyperbolic projective lifts of $P\Gamma$ to theta characteristics. Since the two sets, the set of hyperbolic projective lifts of $P\Gamma$ and the set of theta characteristics, are finite sets with the same cardinality $2^{g}$, it gives rise to a one-to-one correspondence. From the description of the inverse correspondence, the rest of the Theorem follows immediately.
\end{proof}
Note that, for a theta characteristic $\nu$, there is a line bundle $L$ on $X(\Gamma)$ such that $L^{\otimes2}\cong\cO_{X(\Gamma)}$ and $\nu=\omega\otimes L$. We introduce the following (non-standard) definitions which will be used throughout the paper.
\begin{defn2}
On a scheme $X$, a \emph{$2$-torsion line bundle} is a line bundle $L$ over $X$ equipped with an isomorphism $i_{L}:L^{\otimes2}\riso\cO_{X}$ (the isomorphism is a part of the data). Two $2$-torsion line bundles $L,L'$ are regarded as being isomorphic if there exists an isomorphism of line bundles $\iota:L\riso L'$ which is compatible with the isomorphisms $i_{L}:L^{\otimes2}\riso \cO_{X}$ and $i_{L'}:{L'}^{\otimes2}\riso\cO_{X}$. A \emph{trivial $2$-torsion line bundle} is $\cO_{X}$ together with the identity map of $\cO_X$.\footnote{Note that, given a $2$-torsion line bundle $L$ with $i_{L}:L^{\otimes 2}\riso\cO_{X}$, one can scale $i_{L}$ by an invertible element $a\in H^{0}(X,\cO_{X})^{\times}$, and the new $2$-torsion line bundle $(L,ai_{L})$ is isomorphic to $(L,i_{L})$ if and only if $a$ is a square. This is always the case for example if $X$ is a projective variety over an algebraically closed field.}

Similarly, over a field $F$, a \emph{theta characteristic with Kodaira--Spencer data} on $X(\Gamma)_{F}$ is a theta characteristic $\nu$ on $X(\Gamma)_{F}$ equipped with a \emph{Kodaira--Spencer isomorphism} $\KS_{\nu}:\nu^{\otimes 2}\riso\Omega^{1}_{X(\Gamma)_{F}/F}(D)$. Two theta characteristics with Kodaira--Spencer data are isomorphic if there is an isomorphism between the underlying theta characteristics which respects the Kodaira--Spencer isomorphisms.
\end{defn2}
It is well-known that (e.g. \cite[Exercise IV.2.7]{Hartshorne}) there is a one-to-one correspondence between $2$-torsion line bundles and \'etale double covers. Although such a correspondence is usually stated for smooth projective curves over an algebraically closed field, it holds true in much greater generality \emph{if one keeps track of the relevant isomorphisms}.
\begin{lem2}\label{DoubleCoverLem}Let $Y$ be a connected scheme over which $2$ is invertible. Then, there is a natural one-to-one correspondence between the isomorphism classes of $2$-torsion line bundles on $X$ and the isomorphism classes of Galois\footnote{In this paper, by a Galois cover we mean a finite \'etale cover whose self-product splits as a trivial cover. In particular, we allow disconnected covers to be considered Galois.} covers $f:X\rar Y$ of degree $2$. Furthermore, the one-to-one correspondence is compatible with the base-change of $Y$ on both sides.
\end{lem2}
\begin{proof}
Given a $2$-torsion line bundle $L$ over $Y$ (with an isomorphism $i_{L}:L^{\otimes2}\riso\cO_{Y}$), one can define a finite $\cO_{Y}$-algebra \[\cA_{L}:=\cO_{Y}\oplus L,\]such that the multiplication $\cA_{L}\otimes_{\cO_{Y}}\cA_{L}\rar\cA_{L}$ is given by the obvious rules and the isomorphism $i_{L}:L\otimes_{\cO_{Y}}L\rar\cO_{Y}$. Let $X_{L}:=\RelSpec_{\cO_{Y}}\cA_{L}$, which by construction is a finite surjective $Y$-scheme of degree $2$. To show that the natural morphism $X_{L}\rar Y$ is \'etale, we can reduce to the case when $Y=\Spec R$ is affine and $L=\cO_{Y}$. Then, $X_{L}=\Spec R[t]/(t^{2}-a)$ for some $a\in R^{\times}$, which is clearly \'etale over $Y=\Spec R$ as $2t$ is invertible in $R[t]/(t^{2}-a)$. The construction $(L,i_{L})\mapsto X_{L}$ is well-defined, as isomorphic $2$-torsion line bundles yield isomorphic $\cO_{Y}$-algebras. Finally, a morphism $c:X_{L}\rar X_{L}$ induced by the morphism of $\cO_{Y}$-algebras,
\[\cO_{Y}\oplus L\xrar{(x,y)\mapsto (x,-y)}\cO_{Y}\oplus L,\]gies a nontrivial element of $\Aut_{Y}(X_{L})$, which shows that $X\rar Y$ is Galois.

Conversely, suppose that $f:X\rar Y$ is a finite \'etale cover of degree $2$. We claim that $f$ is a Galois covering, i.e. there is a nontrivial involution $\sigma:X\rar X$ of $Y$-schemes. If $X$ is not connected, then it is clear due to the degree reasons that $X$ is consisted of two connected components, each isomorphic to $Y$, so there is an obvious involution of $X$. Thus, suppose that $X$ is connected. Then, the projection to the first coordinate $\pr_{1}:X\times_{Y}X\rar X$ is also a finite \'etale cover of degree $2$ that has a section given by the diagonal $X\rar X\times_{Y}X$. This implies that $X\times_{Y}X$ has a connected component that is isomorphic to $X$ (e.g. \cite[Proposition 5.3.1]{Szamuely}) as an $X$-scheme. Again, by the degree reasons, it follows that $X\times_{Y}X\cong X\coprod X$ as $X$-schemes. Therefore, there is a nontrivial involution $\tau:X\times_{Y}X\rar X\times_{Y}X$ of $X$-schemes (exchanging connected components). By the faithfully flat descent, $\pr_{2}\circ\tau:X\times_{Y}X\rar X$ factors through $\pr_{2}:X\times_{Y}X\rar X$ via a map $\sigma:X\rar X$ that is an involution of $Y$-schemes, as desired.

The involution $\sigma$ induces an endomorphism of $f_{*}\cO_{X}$, which is a vector bundle of rank $2$ over $Y$, and as $2$ is invertible on $Y$, $f_{*}\cO_{X}=(f_{*}\cO_{X})^{\sigma=1}\oplus(f_{*}\cO_{X})^{\sigma=-1}$ where  $(f_{*}\cO_{X})^{\sigma=\pm1}$ is also a vector bundle. Using the descent along $\pr_{2}:X\times_{Y}X\rar X$, we deduce that $(f_{*}\cO_{X})^{\sigma=-1}=:L_{X}$ is a line bundle over $Y$, and $(f_{*}\cO_{X})^{\sigma=1}=\cO_{Y}$; the latter identification is canonical via the adjunction morphism $\cO_{Y}\rar f_{*}\cO_{X}$. 

Note also that the involution $\sigma$ respects the multiplication morphism $(f_{*}\cO_{X})\otimes_{\cO_{Y}}(f_{*}\cO_{X})\rar f_{*}\cO_{X}$, so that the image of its restriction to $L_{X}\otimes_{\cO_{Y}} L_{X}$ lies in $(f_{*}\cO_{X})^{\sigma=1}=\cO_{Y}$. This restriction morphism $L_{X}\otimes_{\cO_Y}L_{X}\rar \cO_Y$ is an isomorphism as it is an isomorphism after base-changing along $\pr_{2}:X\times_{Y}X\rar X$. Thus, given a finite \'etale cover $f:X\rar Y$ of degree $2$, one obtains a $2$-torsion line bundle $L_{X}$ and $L_{X}\otimes_{\cO_Y}L_{X}\riso \cO_Y$. It is clear that the two above constructions are inverses to each other, which implies that these establish a one-to-one correspondence between $2$-torsion line bundles and finite \'etale covers of degree $2$. It is also clear from the constructions that the correspondence is compatible with the base-change of $Y$.
\end{proof}
The correspondence in Lemma \ref{DoubleCoverLem} gives another geometric way to compute $\Gamma_{\nu}$.
\begin{prop2}\label{DoublecoverC}For a $2$-torsion line bundle $L$ on $X(\Gamma)$, let $\alpha:\wt{X}_{L}\rar X(\Gamma)$ be the \'etale double cover corresponding to $L$ under the correspondence of Lemma \ref{DoubleCoverLem}. If we define the representation $\rho_{L}$ of $P\Gamma=\pi_{1}(Y(\Gamma),*)$ to be the composition\[\pi_{1}(Y(\Gamma),*)\rar\pi_{1}(X(\Gamma),*)\thrar\Gal(\wt{X}_{L}/X(\Gamma))=\lbrace\pm1\rbrace,\]the local system $\rho_{\Gamma_{\nu}}$ satisfies $\rho_{\Gamma_{\nu}}=\rho_{\Gamma_{\omega}}\otimes\rho_{L}$. In particular,\[\Gamma_{\nu}=\langle\epsilon_{1}^{a}A_{1},\epsilon_{1}^{b}B_{1},\epsilon_{2}^{a}A_{2},\epsilon_{2}^{b}B_{2},\cdots,\epsilon_{g}^{a}A_{g},\epsilon_{g}^{b}B_{g},C_{1},\cdots,C_{n}\rangle,\]where $\epsilon_{i}^{a}=\rho_{L}(a_{i})$ and $\epsilon_{i}^{b}=\rho_{L}(b_{i})$.
\end{prop2}
Note that, in the above case, the choice of an isomorphism $L^{\otimes 2}\riso \cO_{X(\Gamma)}$ is irrelevant, as $H^{0}(X(\Gamma),\cO_{X(\Gamma)})=\bC$ is an algebraically closed field.
\begin{proof}
As first noted by Deligne, the nonabelian Hodge correspondence is compatible with tensor products (see \cite[p. 8]{SimpsonTensor}). Thus, we only need to show that $\rho_{L}:\pi_{1}(X(\Gamma),*)\rar\lbrace\pm1\rbrace$ and the Higgs bundle $(L,0)$ with zero Higgs field on $X(\Gamma)$ correspond to each other via the nonabelian Hodge correspondence on $X(\Gamma)$. 

Let $c\in\Gal(\wt{X}_{L}/X(\Gamma))$ be the nontrivial Galois element, which gives rise to an automorphism $c\in\Aut_{X(\Gamma)}(\wt{X}_{L})$. Consider $\sH^{0}_{\dR}(\wt{X}_{L}/X(\Gamma))$, which is a vector bundle with an integrable connection $\nabla_{\GM}$ of rank $2$ on $X(\Gamma)$. 
 It is isomorphic to
\[(\sH^{0}_{\dR}(\wt{X}_{L}/X(\Gamma)),\nabla_{\GM})\cong(\cO_{X(\Gamma)},d)\oplus(L,d),\]where $(\cO_{X(\Gamma)},d)$ denotes the canonical differential $d:\cO_{X(\Gamma)}\rar\Omega^{1}_{X(\Gamma)/\bC}$, and $(L,d)=L\otimes(\cO_{X(\Gamma)},d)$ (this defines a connection as $L^{\otimes 2}\cong\cO$, so the transition functions for a sufficiently fine atlas can be taken to be constant functions, namely $\pm1$). Furthermore, $c$ gives rise to an endomorphism of $\sH^{0}_{\dR}(\wt{X}_{L}/X(\Gamma))$, where 
\[\left(\sH^{0}_{\dR}(\wt{X}_{L}/X(\Gamma)),\nabla_{\GM}\right)^{c=1}=(\cO_{X(\Gamma)},d),\quad\left(\sH^{0}_{\dR}(\wt{X}_{L}/X(\Gamma)),\nabla_{\GM}\right)^{c=-1}=(L,d).\]Thus, $\rho_{L}$ (considered as a character) is a local system that underlies a variation of Hodge structure corresponding to $\left(\sH^{0}(\wt{X}_{L}/X(\Gamma))\right)^{c=-1}$, and its associated graded is $(L,0)$. This implies that $\rho_{L}$ and $(L,0)$ correspond to each other via the nonabelian Hodge correspondence.
\end{proof}
We will see later that without much difficulty the same construction works motivically.
\section{The Hodge bundle is the unique congruence theta characteristic}
By Theorem \ref{MainTheta}, for each theta characteristic $\nu$, $H^{0}(X(\Gamma),\nu^{\otimes k})$ is a space of weight $k$ modular forms of level $\Gamma_{\nu}$. Note that, if $k$ is even, then \[\nu^{\otimes k}=\left(\nu^{\otimes 2}\right)^{\otimes k/2}=\left(\omega^{\otimes 2}\right)^{\otimes k/2}=\omega^{\otimes k},\]so $M_{k}(\Gamma_{\nu})=M_{k}(\Gamma)$ is consisted of modular forms of level $\Gamma$ which is a congruence subgroup. 

From this, one naturally wonders about the nature of \emph{odd-weight} modular forms of level $\Gamma_{\nu}$. We introduce the following definition.
\begin{defn2}
A theta characteristic $\nu$ is called a \emph{congruence theta characteristic} if $\Gamma_{\nu}\le\SL_{2}(\bZ)$ is a congruence subgroup.
\end{defn2}
In contrast to the even-weight case, the main theorem of this section shows that, if $k$ is odd and $\nu\ne\omega$, $M_{k}(\Gamma_{\nu})$ is consisted entirely of \emph{noncongruence modular forms}\footnote{Recall that a \emph{noncongruence modular form} is a modular form of some level which does not arise as a modular form of congruence level.}!
\begin{thm2}\label{CongruenceTheta}
A theta characteristic $\nu$ is a congruence theta characteristic if and only if $\nu=\omega$.
\end{thm2}
A quick corollary is that, for $\nu\ne\omega$, the Hecke operators are zero on $H^{0}(X(\Gamma),\nu^{\otimes k})$ for odd $k$.
\begin{coro2}\label{HeckeZero}
For $(p,N)=1$ and odd $k\ge1$, define the Hecke operator $T_{p}$ on $H^{0}(X(\Gamma),\nu^{\otimes k})=M_{k}(\Gamma_{\nu})$ as follows: for $f\in M_{k}(\Gamma_{\nu})$,
\[T_{p}f=\sum_{i}f\rvert_{\alpha\alpha_{i}},\]where \[\Gamma_{\nu}\alpha\Gamma_{\nu}=\cup_{i}\Gamma_{\nu}\alpha\alpha_{i},\quad \alpha=\begin{pmat}p&0\\0&1\end{pmat}.\]If $\nu\ne\omega$, we always have $T_{p}f=0$.
\end{coro2}
\begin{proof}
By \cite{Berger}, we know that $T_{p}$ factors through the trace map to the congruence closure. In our case, if $\nu\ne\omega$, by Theorem \ref{CongruenceTheta}, the congruence closure of $\Gamma_{\nu}$ is $\langle\pm1,\Gamma\rangle$. Since there is no nonzero odd-weight modular form of level $\langle\pm1,\Gamma\rangle$, the desired  statement follows.
\end{proof}
\begin{rmk2}
The above Hecke operator can be geometrically interpreted as the correspondence
\[\xymatrix{&X(\Gamma_{\nu}\cap\alpha^{-1}\Gamma_{\nu}\alpha\cap\Gamma_{0}(p))\ar[ld]\ar[r]^-{\sim}&X(\alpha\Gamma_{\nu}\alpha^{-1}\cap\Gamma_{\nu}\cap\Gamma^{0}(p))\ar[rd]&\\X(\Gamma_{\nu})&&&X(\Gamma_{\nu})}\]where $\Gamma^{0}(p)=\left\lbrace\begin{psmat}a&b\\c&d\end{psmat}\equiv\begin{psmat}*&0 \\ *&*\end{psmat}~(\MOD p)\right\rbrace$. As $\Gamma_{\nu}$ is an index $2$ subgroup of a congruence subgroup of level prime to $p$, $\alpha$ in general does not normalize $\Gamma_{\nu}$, but rather sends $\Gamma_{\nu}$ to a possibly different index $2$ subgroup of $\langle\pm1,\Gamma\rangle$.
\end{rmk2}
The proof of Theorem \ref{CongruenceTheta} will be a slight generalization of the proofs in \cite[\S2]{Kiming2}, for which we use the technical condition ($*$). As in \emph{loc. cit.}, we define \[V_{2}(G):=G^{\ab}/(G^{\ab})^{2}=G^{\ab}\otimes_{\bZ}\bF_{2},\]for a group $G$. Note that $V_{2}$ is a functor that sends finitely generated groups to finite-dimensional $\bF_{2}$-vector spaces, which has the following easy properties.
\begin{lem2}
\hfill
\begin{enumerate}
\item $V_{2}$ is a right-exact functor.\footnote{Even though the category of groups is not an abelian category, the notion of exact sequences makes sense.}
\item $V_{2}(G_{1}\times G_{2})\cong V_{2}(G_{1})\times V_{2}(G_{2})$.
\end{enumerate}
\end{lem2}
\begin{proof}
The functor $V_{2}$ is the composition of the abelianization functor with $(-)\otimes_{\bZ}\bF_{2}$, and both are right exact.
\end{proof}
\begin{proof}[Proof of Theorem \ref{CongruenceTheta}]
As $\Gamma_{\omega}=\Gamma$, $\Gamma_{\omega}$ is a congruence subgroup, which proves one direction. Conversely, suppose that $\Gamma_{\nu}$ is a congruence subgroup. As per Theorem \ref{MainTheta}, we need to prove that there is no hyperbolic projective lift of $P\Gamma$ which is a congruence subgroup and is different from $\Gamma$. 

Suppose that $\Gamma$ is of level $N$; namely, $N$ is the minimal number such that $\Gamma(N)\le\Gamma$. By the result of Wohlfahrt \cite[Theorem 2]{Wohlfahrt} and Kiming--Sch\"utt--Verrill \cite[Proposition 3]{Kiming1}, $\Gamma$ is of general level either $N$ or $\frac{N}{2}$. Recall that the \emph{general level} of a Fuchsian group is the least common multiple of the widths of the cusps. The general level only depends on the projective image of the Fuchsian group, so $\Gamma_{\nu}$ is of general level $N$. By \emph{loc. cit.}, $\Gamma_{\nu}\ge\Gamma(2N)$. Thus, $\langle\pm1,\Gamma\rangle\ge\Gamma_{\nu}\ge\Gamma(2N)$. Thus, $\Gamma_{\nu}$ corresponds to a subgroup of $\langle\pm1,\Gamma\rangle/\Gamma(2N)\cong\langle\pm1\rangle\times\Gamma/\Gamma(2N)$ such that $\lbrace\pm1\rbrace$ and $\Gamma_{\nu}$ together generate the whole subgroup $\langle\pm1,\Gamma\rangle/\Gamma(2N)$. As in \cite[Proposition 1]{Kiming2}, projective lifts of $P\Gamma$ that are also congruence subgroups are in one-to-one corresopndence with a sub-$\bF_{2}$-vector space $U$ of $V_{2}(\langle\pm1,\Gamma\rangle/\Gamma(2N))\cong\langle\pm1\rangle\times V_{2}(\Gamma/\Gamma(2N))$ such that $U$ and $-1$ together span the whole vector space. Such projective lift is honest if $U$ is a proper subspace, and $-1\notin U$. Thus, the composition $U\hrar V_{2}(\langle\pm1,\Gamma\rangle/\Gamma(2N))=\langle\pm1\rangle\times V_{2}(\Gamma/\Gamma(2N))\thrar V_{2}(\Gamma/\Gamma(2N))$ is injective, thus bijective (as the target and the source have the same $\bF_{2}$-dimensions). Thus, choosing an honest congruence projective lift is the same as choosing the signs for the lifts of basis elements of $V_{2}(\Gamma/\Gamma(2N))$.

By the assumption ($*$), $N=\lcm(N_{1},N_{2})$, and $\Gamma(N)\le\Gamma\le\Gamma_{1}(N)$. Let $N=2^{s}p_{1}^{t_{1}}\cdots p_{r}^{t_{r}}$, where $p_{i}$'s are odd primes. Note also that
\[\SL_{2}(\bZ)/\Gamma(2N)\cong\SL_{2}(\bZ/2N\bZ)\cong\SL_{2}(\bZ/2^{s+1}\bZ)\times\prod_{i=1}^{r}\SL_{2}(\bZ/p_{i}^{t_{i}}\bZ),\]so $\Gamma_{1}(N)/\Gamma(2N)$ injects into $\Gamma_{1}(2^{s})/\Gamma(2^{s+1})\times\prod_{i=1}^{r}\Gamma_{1}(p_{i}^{t_{i}})/\Gamma(p_{i}^{t_{i}})$, which is a bijection as the two groups are finite groups of the same order; for $1\le a\le b$, $\#\Gamma_{1}(p^{a})/\Gamma(p^{b})=p^{3b-2a}$. 
Under this isomorphism, we have
\[\Gamma/\Gamma(2N)\cong A\times\prod_{i=1}^{r}B_{i},\]where $B_{i}\le\Gamma_{1}(p_{i}^{t_{i}})/\Gamma(p_{i}^{t_{i}})$ is a subgroup, and \[A=\begin{cases}\SL_{2}(\bZ)/\Gamma(2)&\text{ if $s=0$,}\\\Gamma_{1}(2^{s})/\Gamma(2^{s+1})&\text{ if $2|N_{1}$,}\\\Gamma(2^{s})/\Gamma(2^{s+1})&\text{ if $2|N_{2}$.}\end{cases}\]Note that $\Gamma_{1}(p_{i}^{t_{i}})/\Gamma(p_{i}^{t_{i}})$ is of odd order, so $B$ is of odd order as well. Thus, the natural projection map $\Gamma/\Gamma(2N)\thrar A$ induces an isomorphism $V_{2}(\Gamma/\Gamma(2N))\riso V_{2}(A)$. 

By the right-exactness of $V_{2}$, we have a surjective natural map $V_{2}(\Gamma)\thrar V_{2}(\Gamma/\Gamma(2N))$. Since a hyperbolic projective lift fixes the signs of the lifts of the loops around the cusps, to prove Theorem \ref{CongruenceTheta}, it suffices to prove that $V_{2}(\Gamma/\Gamma(2N))$ is spanned by the images of shearing transformations along the cusps. We prove that this is true by dividing into cases.
\begin{enumerate}
\item[(Case 1)] If $s=0$, then $A=\SL_{2}(\bF_{2})\cong S_{3}$, and $V_{2}(A)\cong(\bZ/2\bZ)$ is generated by $\begin{psmat}1&1\\0&1\end{psmat}$. As $\Gamma=\Gamma_{1}(N_{1})\cap\Gamma(N_{2})$ with $N_{1},N_{2}$ odd, $\begin{psmat}1&N_{2}\\0&1\end{psmat}\in\Gamma$ is sent to $\begin{psmat}1&1\\0&1\end{psmat}\in A$ via the natural projection $\Gamma\thrar A$. Since $\begin{psmat}1&N_{2}\\0&1\end{psmat}$ is a shearing transformation along the cusp $\infty\in\bP^{1}(\bQ)$, there is no hyperbolic projective lift different from $\Gamma$.
\item[(Case 2)] If $s>0$ and $2|N_{2}$, then $A=\Gamma(2^{s})/\Gamma(2^{s+1})$. As in the proof of \cite[Proposition 2]{Kiming2}, one notes that $A=V_{2}(A)\cong(\bZ/2\bZ)^{3}$ with a generator given by
\[\alpha=\begin{pmat}1&2^{s}\\0&1\end{pmat},\quad \beta=\begin{pmat}1+2^{s}&-2^{s}\\2^{s}&1-2^{s}\end{pmat}\quad\gamma=\begin{pmat}1&0\\2^{s}&1\end{pmat}.\]Note that we took a slightly different set of generators. Note that
\[\alpha\equiv\begin{pmat}1&N_{2}\\0&1\end{pmat}(\MOD 2^{s+1}),\quad \beta\equiv\begin{pmat}1+N_{1}N_{2}&-N_{1}N_{2}\\N_{1}N_{2}&1-N_{1}N_{2}\end{pmat}(\MOD 2^{s+1}),\quad \gamma\equiv\begin{pmat}1&0\\N_{1}N_{2}&1\end{pmat}(\MOD 2^{s+1}),\]and these matrices are genuine elements of $\Gamma=\Gamma_{1}(N_{1})\cap\Gamma(N_{2})$. Also note that $\begin{psmat}1&N_{2}\\0&1\end{psmat}$ is a shearing transformation along the cusp $\infty\in\bP^{1}(\bQ)$, $\begin{psmat}1&0\\N_{1}N_{2}&1\end{psmat}$ is a shearing transformation along the cusp $0\in\bP^{1}(\bQ)$, and $\begin{psmat}1+N_{1}N_{2}&-N_{1}N_{2}\\N_{1}N_{2}&1-N_{1}N_{2}\end{psmat}$ is a shearing transformation along the cusp $1\in\bP^{1}(\bQ)$, since
\[\begin{pmat}1+N_{1}N_{2}&-N_{1}N_{2}\\N_{1}N_{2}&1-N_{1}N_{2}\end{pmat}=\begin{pmat}1&1\\0&1\end{pmat}\begin{pmat}1&0\\N_{1}N_{2}\end{pmat}\begin{pmat}1&-1\\0&1\end{pmat}.\]Therefore, there is no hyperbolic projective lift different from $\Gamma$.
\item[(Case 3)] If $s>0$ and $2|N_{1}$, then $A=\Gamma_{1}(2^{s})/\Gamma(2^{s+1})$. As per \emph{loc. cit.}, $V_{2}(A)\cong(\bZ/2\bZ)^{2}$ with a basis given by \[\tau=\begin{pmat}1&1\\0&1\end{pmat},\quad\gamma=\begin{pmat}1&0\\2^{s}&1\end{pmat}.\]Since $N_{2}$ is odd, $2^{s+1}$ is invertible modulo $N_{2}$, which implies that there exists $k\in\bZ$ such that $k2^{s+1}\equiv-1(\MOD N_{2})$. Now note that 
\[\tau\equiv\begin{pmat}1&1+k2^{s+1}\\0&1\end{pmat}(\MOD 2^{s+1}), \quad \gamma\equiv\begin{pmat}1&0\\N_{1}N_{2}&1\end{pmat}(\MOD 2^{s+1}),\]and these matrices are genuine elements of $\Gamma=\Gamma_{1}(N_{1})\cap\Gamma(N_{2})$. It is clear that these matrices are also shearing transformations along the cusps $\infty,0\in\bP^{1}(\bQ)$, respectively, so there is no hyperbolic projective lift different from $\Gamma$.
\end{enumerate}
\end{proof}
\section{Brill--Noether theory of the modular curves and the Hodge bundle}
There is a different perspective on how special a line bundle on a curve is, which goes under the general name of \emph{Brill--Noether theory}. Generally speaking, given a smooth projective complex curve $C$ and a line bundle $L$ (or equivalently a divisor class $[D]$), $L$ is considered  special if $\dim_{\bC}H^{0}(C,L)$ is larger than the other line bundles on $C$ of the same degree. Furthermore, the curve $C$ is considered special if there exists a certain line bundle with a larger than usual $h^{0}$. More precisely, we introduce the following definition which is common in the literature.
\begin{defn2}\label{BNGeneral}
A smooth projective complex curve $C$ of genus $g$ is called \emph{Brill--Noether general} if, for all $r,d\ge0$, the moduli space $G_{d}^{r}(C)$ of linear systems on $C$ of degree $d$ and dimension $r$ (for the definition, see \cite[pg. 177]{ACGH}) is smooth of dimension $\rho(g,d,r):=g-(r+1)(g-d+r)$; this in particular means that $G_{d}^{r}(C)$ is empty if $\rho(g,d,r)<0$. 

Otherwise, $C$ is called \emph{Brill--Noether special}.
\end{defn2}
The classical Brill--Noether theory (\cite{Kempf}, \cite{KleimanLaksov}, \cite{GriffithsHarris}, \cite{GiesekerPetri}) shows that a general curve is Brill--Noether general (i.e. there is a Zariski dense open subset of the moduli $\cM_{g}$ of genus $g$ curves whose closed points are Brill--Noether general curves). Given the special nature of the modular curves and their Hodge bundles, one is naturally led to ask whether they are special in the sense of Brill--Noether theory. 

We first show that a modular curve of sufficiently fine level is Brill--Noether special. For this, we notice a simple lemma, which extends the idea of \cite[Remark 1.12]{Welters}.
\begin{lem2}\label{BNspecial}
Let $C$ be a smooth projective complex curve, and suppose that there exists a degree $d$ divisor $D$ satisfying the following two conditions:
\begin{enumerate}
\item $h^{0}(\cO(D))\ge2$;
\item $K_{C}-2D$ is linearly equivalent to an effective divisor (may be zero).
\end{enumerate}
If $\wt{C}\rar C$ is a finite surjective morphism from a smooth projective complex curve $\wt{C}$, then $\wt{C}$ is Brill--Noether special.
\end{lem2}
\begin{proof}
We first show that $C$ itself is Brill--Noether special, and then show that $\wt{C}$ itself has a divisor satisfying the two conditions. Take a $2$-dimensional subspace $V\subset H^{0}(C,\cO(D))$. We will show that $G_{d}^{1}(C)$ is not smooth at the point corresponding to the linear system $(D,V)$, which will then show that $C$ is Brill--Noether special. This will follow if the Petri map
\[\mu_{D,V}:V\otimes H^{0}(C,\Omega_{C}^{1}(-D))\rar H^{0}(C,\Omega_{C}^{1}),\]which is simply the cup product, is not injective. Note that $\Omega_{C}^{1}(-D)\cong \cO(D+E)$, where $E$ is an effective divisor. Therefore, any section $s\in H^{0}(C,\cO(D))$ is also a section of $\cO(D+E)$. Take two linearly independent sections $s,s'\in V$, then $s\otimes s'-s'\otimes s$ is in $\ker\mu_{D,V}$, which implies that the above Petri map is indeed not injective, as desired. 

Suppose now that $\pi:\wt{C}\rar C$ is a finite surjective morphism. We claim that $\pi^{*}D$ satisfies the properties (1) and (2). By \cite[(2.2.8)]{EGA4b}, $H^{0}(C,\cO(D))\rar H^{0}(\wt{C},\cO(\pi^{*}D))$ is injective, so (1) holds. If $K_{C}\sim 2D+E$ for some effective divisor $E$, then $\pi^{*}K_{C}\sim 2\pi^{*}D+\pi^{*}E$. By Hurwitz's theorem (e.g. \cite[Proposition IV.2.3]{Hartshorne}), $K_{\wt{C}}\sim 2\pi^{*}D+(\pi^{*}E+R)$ where $R$ is the ramification divisor, which is in particular effective. Thus, (2) holds as well.
\end{proof}
\begin{prop2}\label{ModCurveBN}
The modular curve $X(\Gamma)$ is Brill--Noether special if one of the following conditions holds:
\begin{enumerate}
\item $\dim_{\bC} S_{1}(\Gamma)\ge2$;
\item $\Gamma=\Gamma_{1}(N)$, when $N$ is a squarefree number such that $N>10$ and $4|\varphi(N)$;
\item $\Gamma\le\Gamma'$ for $\Gamma'$ satisfying any of the above conditions.
\end{enumerate}
\end{prop2}  
\begin{proof}
\begin{enumerate}
\item The Petri map for $\omega(-D)$ is \[\mu_{\omega(-D)}:H^{0}(X(\Gamma),\omega(-D))\otimes H^{0}(X(\Gamma),\omega)\rar H^{0}(X(\Gamma),\Omega_{C}^{1}),\]which is simply the multiplication map $S_{1}(\Gamma)\otimes M_{1}(\Gamma)\rar S_{2}(\Gamma)$. For any two-dimensional subspace $V\subset S_{1}(\Gamma)$, $s\otimes s'-s'\otimes s$ is in the kernel of the Petri map, where $s,s'$ are two linearly independent sections of $V$, which implies that $G_{g-1-\frac{n}{2}}^{1}(X(\Gamma))$ is not smooth at the point corresponding to the linear system $(\omega(-D),V)$.
\item The map $X_{1}(N)\rar X_{0}(N)$ is a Galois cover with the Galois group $(\bZ/N\bZ)^{\times}/\lbrace\pm1\rbrace$. By the assumption on $N$, the Galois group is of even order, so one may find an order $2$ element in the Galois group, which corresponds to an \'etale double cover $X_{1}(N)\rar X$ for some smooth projective complex curve $X$ of genus $>0$. By \cite[Remark 1.12]{Welters}, $X_{1}(N)$ is Brill--Noether special.  
\item In the above two cases, $X(\Gamma')$ is shown to be Brill--Noether special using a divisor satisfying the conditions (1) and (2) of Lemma \ref{BNspecial}. Thus, we can apply  Lemma \ref{BNspecial} to the covering $X(\Gamma)\rar X(\Gamma')$ to deduce that $X(\Gamma)$ is also Brill--Noether special.
\end{enumerate}
\end{proof}
In contrast to the modular curves, it turns out that the Hodge bundle $\omega$ has no meaningful relationship with the Brill--Noether theory, as illustrated by the following two examples.
\begin{exam2}[(When $\omega$ is the theta characteristic with the most sections)]\label{SpecialW}
Let $\Gamma=\Gamma_{1}(23)$. According to \cite{LMFDB}, $\dim_{\bC} H^{0}(X_{1}(23),\omega(-D))=\dim_{\bC}S_{1}(\Gamma_{1}(23))=1$. On the other hand, $X_{1}(23)$ is of genus $12$ and has $22$ cusps, so $\deg\omega(-D)=0$. Therefore, $\omega(-D)\cong\cO_{X_{1}(23)}$. For any other theta characteristic $\nu$, $\nu(-D)$ is of degree zero and is not the structure sheaf, so $H^{0}(X_{1}(23),\nu(-D))=0$. Therefore, $\dim_{\bC}H^{0}(X_{1}(23),\omega)>\dim_{\bC}H^{0}(X_{1}(23),\nu)$ for any theta characteristic $\nu\ne\omega$, or $\omega$ is the ``most special'' theta characteristic.
\end{exam2}
\begin{exam2}[(When $\omega$ is the theta characteristic with the fewest sections)]\label{GeneralW} Let $\Gamma=\Gamma_{1}(35)$. According to \cite{LMFDB}, $\dim_{\bC}H^{0}(X_{1}(35),\omega(-D))=\dim_{\bC}S_{1}(\Gamma_{1}(35))=0$. We claim that there is a unique theta characteristic $\nu$ such that $\dim_{\bC}H^{0}(X_{1}(35),\nu(-D))=1$, and any other theta characteristic $\nu'\ne\nu$ (including $\nu'=\omega$) satisfies $\dim_{\bC}H^{0}(X_{1}(35),\nu'(-D))=0$. This implies that $\omega$ is the ``least special'' theta characteristic in this case.

Note that $X_{1}(35)$ is of genus $25$ and has $48$ cusps, so in particular $\deg\Omega^{1}_{X_{1}(35)}(-D)=0$. This  means that if there exists a nonzero weight $2$ cusp form on $X_{1}(35)$ which vanishes to order $2$ at all cusps, then it does not vanish outside the cusps and therefore is a square of a holomorphic function, which is a weight one non-congruence modular form. Furthermore, if $\nu$ is the theta characteristic corresponding to the level of this non-congruence modular form, then as $\deg\nu(-D)=0$, $\nu(-D)\cong\cO_{X_{1}(35)}$, which implies that $\dim_{\bC}H^{0}(X_{1}(35),\nu(-D))=1$. Moreover, for any other theta characteristic $\nu'\ne\nu$, $\nu'(-D)$ is a degree zero line bundle that is not the structure sheaf, which will imply that $\dim_{\bC}H^{0}(X_{1}(35),\nu'(-D))=0$.

Let $f_{1},f_{2}\in S_{2}(\Gamma_{1}(35))$ be the embedded newforms 35.2.a.b.1.1 and 35.2.a.b.1.2 in \cite{LMFDB}, respectively. We claim that
$f:=f_{2}-f_{1}$ vanishes to order $2$ at all cusps of $X_{1}(35)$. We may compute its $q$-expansion at the  cusps other than the infinity cusp using \cite{Asai}. To state the results, we introduce some notations. For the remainder of this example, we assume that $N=p_{1}\cdots p_{r}$ is a square-free odd integer. 

There are $2^{r-1}(p_{1}-1)\cdots(p_{r}-1)$ cusps of $X_{1}(N)$, where we can take the representatives of the cusps nicely as
\[\left\lbrace c_{M_{1},a,b}:=\frac{b}{M_{2}a}\in\bP^{1}(\bQ)~:~N=M_{1}M_{2},~1\le a<M_{1},~(a,M_{1})=1,~(b,M_{2}a)=1\right\rbrace/\sim,\]where $c_{M_{1},a,b}\sim c_{M_{1}',a',b'}$ if and only if $M_{1}=M_{1}'$ and if there exists $\epsilon\in\lbrace\pm1\rbrace$ such that $a\equiv \epsilon a'~(\MOD M_{1})$ and $b\equiv \epsilon b'~(\MOD M_{2})$. Note that $c_{1,1,1}$ is the infinity cusp, and $c_{M_{1},a,b}$ is of width $M_{1}$. For a modular form $g\in M_{k}(\Gamma_{1}(N))$, let
\[e_{M_{1},a,b}(g)(q)=\sum_{j=0}^{\infty}a_{M_{1},a,b,j}(g)q^{j/M_{1}}\in\bC[[q^{1/M_{1}}]],\]be the $q$-expansion of $g$ at the cusp $c_{M_{1},a,b}$; more precisely, it is the Fourier expansion of $f_{M_{1},a,b}(z):=f(\sigma_{M_{1},a,b}(z))$, where $\sigma_{M_{1},a,b}\in\SL_{2}(\bZ)$ is a matrix such that \[\sigma_{M_{1},a,b}^{-1}Z_{c_{M_{1},a,b}}(\Gamma_{1}(N)),\sigma_{M_{1},a,b}\subset\begin{pmat}1& *\\0&1\end{pmat},\]where $Z_{c_{M_{1},a,b}}(\Gamma_{1}(N))\subset \Gamma_{1}(N)$ is the stabilizer of $c_{M_{1},a,b}$ in $\Gamma_{1}(N)$.

We may now state a result of Asai computing the $q$-expansion of a new Hecke eigenform at all cusps.
\begin{prop2}[({\cite[Theorem 2]{Asai}})]\label{Newform}
Let $g\in S_{k}^{\new}(N,\chi)$ be a Hecke eigenform such that, at $c_{1,1,1}=\infty$, $g$ has the $q$-expansion\[e_{1,1,1}(g)(q)=\sum_{n=1}^{\infty}a_{n}q^{n}, \qquad a_{1}=1.\]Let $\chi=\chi_{p_{1}}\cdots\chi_{p_{r}}$ be the product such that $\chi_{p_{i}}$ is a Dirichlet character mod $p_{i}$ (may not be primitive). Then, for all cusps $c_{M_{1},a,b}$,\[e_{M_{1},a,b}(g)(q)=\chi\left(bcM_{1}+M_{2}^{2}ad\right)\prod_{p|M_{1}}\left(p^{-\frac{k}{2}}\chi_{p}\left(\frac{M_{1}}{p}\right)\ov{a_{p}}C(\chi_{p})\right)\sum_{n=1}^{\infty}a_{n}^{(M_{1})}q^{n/M_{1}},\]where $c,d\in\bZ$ are such that $cM_{1}+dM_{2}=1$, 
\[C(\chi_{p})=\begin{cases}\sum_{1\le h<p}\chi_{p}(h)e^{2\pi ih/p}&\text{ if $\chi_{p}$ is primitive}\\-q&\text{ if $\chi_{p}$ is trivial,}\end{cases}\]and $a_{n}^{(M_{1})}$ is defined as
\[a_{n}^{(M_{1})}=\begin{cases}\ov{\chi(dnM_{2}+cM_{1})}a_{n}&\text{ if $(n,M_{1})=1$}\\\chi(cnM_{1}+dM_{2})\ov{a}_{n}&\text{ if $\left(n,M_{2}\right)=1$}\\a_{x}^{(M_{1})}a_{y}^{(M_{1})}&\text{ if $n=xy$, $(x,y)=1$.}\end{cases}\]
\end{prop2}
For our purpose, this can be packaged more simply as follows.
\begin{coro2}\label{NewformQ}Fix $\chi$ and a cusp $c_{M_{1},a,b}\in X_{1}(N)$ of width $M_{1}$. Then, there exist a constant $\lambda\in\bC$ that depends only on $M_{1},a,b$, and, for each $n\ge1$, a constant $\epsilon_{n}\in\bC$ that depends only on $M_{1},a,b,n$, such that, for any normalized Hecke eigenform $g\in S_{k}^{\new}(N,\chi)$ with the $q$-expansion $e_{1,1,1}(g)(q)=\sum_{n=1}^{\infty}a_{n}q^{n}$ at $\infty$, 
\[a_{M_{1},a,b,n}(g)=\begin{cases}\left(\lambda\epsilon_{n}\prod_{p|M_{1}}\ov{a_{p}}\right)a_{n}&\text{if $(n,M_{1})=1$} \\ \left(\lambda\epsilon_{n}\prod_{p|M_{1}}\ov{a_{p}}\right)\ov{a_{n}}&\text{if $(n,M_{2})=1$} \\ \frac{a_{M_{1},a,b,x}(g)a_{M_{1},a,b,y}(g)}{\lambda\epsilon_{1}\prod_{p|M_{1}}\ov{a_{p}}}&\text{if $n=xy$, $(x,y)=1$.}\end{cases}\]
\end{coro2}
Now we are ready to show that $f$ vanishes to order $2$ at all cusps of $X_{1}(35)$. 
Note first that $f_{1},f_{2},f$ all have trivial character, so the order of vanishing of $f$ is constant along the cusps of $X_{1}(35)$ in a single fiber of the Galois covering $X_{1}(35)\rar X_{0}(35)$, namely the order of vanishing of $f$ is constant along the cusps of the same width. As $a_{1,1,1,1}(f)=0$, we only need to show that $a_{5,1,1,1}(f)=a_{7,1,1,1}(f)=a_{35,1,1,1}(f)=0$. Let $f_{i}=\sum_{j=1}^{\infty}a_{i,j}q^{n}$ be the $q$-expansion of $f_{i}$. As $f_{1}$ and $f_{2}$ are both newforms, we may use Corollary \ref{NewformQ} to compute that, for $M_{1}=5,7,35$,
\[a_{M_{1},1,1,1}(f)=C_{M_{1}}\left(\prod_{p|M_{1}}\ov{a_{2,p}}-\prod_{p|M_{1}}\ov{a_{1,p}}\right),\]for some constant $C_{M_{1}}\in\bC$ that only depends on $M_{1}$. According to \cite{LMFDB}, we have
\[a_{1,5}=1, \quad a_{1,7}=-1,\]
\[a_{2,5}=1,\quad a_{2,7}=-1,\]which implies that $a_{M_{1},1,1,1}(f)=0$ for $M_{1}=5,7,35$. This implies that $f$ vanishes to order $2$ at all cusps of $X_{1}(35)$, as desired.
\end{exam2}
\section{Twisted Kuga--Sato varieties and geometric local systems for theta characteristics}\label{KugaSato}
We now aim to show that the construction of Lemma \ref{DoubleCoverLem} yields, for each theta characteristic with Kodaira--Spencer data $\nu$, a compatible system of  {local systems} over $Y(\Gamma)_{K}$ over an appropriate number field $K$ that \emph{comes from geometry}. Furthermore, we will show that these local systems correspond to the uniformizing logarithmic Higgs bundles $(E_{\nu},\theta_{\nu})$ via complex and $p$-adic nonabelian Hodge correspondneces. We will construct the geometric local systems from the \emph{twisted Kuga--Sato variety}. 
\begin{defn2}
Let $K/\bQ$ be a number field, and let $\alpha:\wt{X}\rar X(\Gamma)_{K}$ be a finite  Galois cover of degree $r$. The \emph{twisted Kuga--Sato variety} $u:\ov{W}_{\wt{X}}\rar X(\Gamma)_{K}$ associated with $\wt{X}$ is defined as the Weil restriction of the pullback $\alpha^{*}\ov{\cE}_{K}$,
\[\ov{W}_{\wt{X}}:=\tR_{\wt{X}/X(\Gamma)_{K}}(\alpha^{*}\ov{\cE}_{K}).\]
The open twisted Kuga--Sato variety $u:W_{\wt{X}}\rar Y(\Gamma)_{K}$ is defined as the open subscheme of $\ov{W}_{\wt{X}}$ lying over $Y(\Gamma)_{K}\subset X(\Gamma)_{K}$.

We will also use the notation $\wt{Y}:=\alpha^{-1}(Y(\Gamma)_{K})$.
\end{defn2}
For the definition of scheme-theoretic Weil restriction of scalars, see \cite[\S7.6]{BLR}.
\begin{exam2}
For a trivial Galois $r$-cover $X(\Gamma)_{\bQ}\coprod\cdots\coprod X(\Gamma)_{\bQ}\rar X(\Gamma)_{\bQ}$ consisted of $r$ copies of $X(\Gamma)_{\bQ}$, the corresponding twisted Kuga--Sato variety is the usual Kuga--Sato variety (before the canonical desingularization), namely the $r$-fold fiber product of $\ov{\cE}_{\bQ}$ over $X(\Gamma)_{\bQ}$.
\end{exam2}
\begin{rmk2}
As in the case of the usual Kuga--Sato variety, the twisted Kuga--Sato variety $\ov{W}_{\wt{X}}$ is in general singular, even though $\alpha^{*}\ov{\cE}_{K}$ itself is a smooth $K$-scheme. On the other hand, the open twisted Kuga--Sato variety $W_{\wt{X}}$ is smooth.
\end{rmk2}
The open twisted Kuga--Sato variety ${W}_{\wt{X}}$ still turns out to be a family of principally polarized abelian $r$-folds over $Y(\Gamma)_{K}$.
\begin{prop2}
The open twisted Kuga--Sato variety $u:W_{\wt{X}}\rar Y(\Gamma)_{K}$ is a family of principally polarized abelian varieties of dimension $r$.
\end{prop2}\begin{proof}
By \cite[Proposition 2]{DN}, the Weil restriction of a principal polarization is a principal polarization. As an elliptic curve is canonically principally polarized, $W_{\wt{X}}$ is a family of principally polarized abelian varieties.
\end{proof}
We can thus think of classifying map to the moduli space of principally polarized abelian varieties of dimension $r$,
\[\pi_{\wt{X}}:Y(\Gamma)_{K}\rar\cA_{r,K},\]and this induces a classifying map to the corresponding coarse moduli scheme,
\[p_{\wt{X}}:Y(\Gamma)_{K}\rar A_{r,K}.\]
\begin{prop2}
Let $\pi_{\diag}:Y(\Gamma)_{K}\rar\cA_{r,K}$ be the classifying map which corresponds to the $r$-th self-product $\cE_{K}\times_{Y(\Gamma)_{K}}\cdots\times_{Y(\Gamma)_{K}}\cE_{K}$ over $Y(\Gamma)_{K}$, and let $p_{\diag}:Y(\Gamma)_{K}\rar A_{r,K}$ be the corresponding map to the coarse moduli scheme. Then, $p_{\diag}=p_{\wt{X}}$. On the other hand, $\pi_{\diag}\ne\pi_{\wt{X}}$ unless $\wt{X}\rar X(\Gamma)_{K}$ is a trivial cover. 
\end{prop2}
\begin{proof}By \cite[Proposition A.5.2]{CGP}, we have
\[W_{\wt{X}}\times_{Y(\Gamma)_{K}}\wt{Y}=\tR_{\wt{Y}\times_{Y(\Gamma)_{K}}\wt{Y}/\wt{Y}}(\alpha^{*}\cE_{K}\times_{\wt{Y}}(\wt{Y}\times_{Y(\Gamma)_{K}}\wt{Y})).\]Since $\wt{Y}\times_{Y(\Gamma)_{K}}\wt{Y}$ is isomorphic to the disjoint union of $r$ copies of $\wt{Y}$, we have
\[W_{\wt{X}}\times_{Y(\Gamma)_{K}}\wt{Y}\cong\tR_{\wt{Y}\coprod\cdots\coprod\wt{Y}/\wt{Y}}\left(\alpha^{*}\cE_{K}\coprod\cdots\coprod \alpha^{*}\cE_{K}\right)\cong \alpha^{*}\cE_{K}\times_{\wt{Y}}\cdots\times_{\wt{Y}}\alpha^{*}\cE_{K}.\]
This implies that
\[W_{\wt{X}}\times_{Y(\Gamma)_{K}}\wt{Y}\cong(\cE_{K}\times_{Y(\Gamma)_{K}}\cdots\times_{Y(\Gamma)_{K}}\cE_{K})\times_{Y(\Gamma)_{K}}\wt{Y}.\]Furthermore, the two abelian schemes are isomorphic as $\wt{Y}$-schemes; even though $\wt{Y}\times_{Y(\Gamma)_{K}}\wt{Y}\cong\wt{Y}\coprod\cdots\coprod\wt{Y}$ is most naturally thought as being indexed by the elements in $\Gal(\wt{Y}/Y(\Gamma)_{K})$, for any $\sigma\in\Gal(\wt{Y}/Y(\Gamma)_{K})$, the $\wt{Y}$-scheme $\wt{Y}\xrar{\sigma}\wt{Y}$ is isomorphic as a $\wt{Y}$-scheme to $\wt{Y}\xrar{\id}\wt{Y}$. Thus, $\pi_{\diag}\circ\alpha=\pi_{\wt{X}}\circ\alpha$, and $p_{\diag}\circ\alpha=p_{\wt{X}}\circ\alpha$. Since $\alpha$ is flat, surjective, and locally of finite presentation, by \cite[Tag 05VM]{Stacks}, $\alpha$ is an epimorphism (i.e. surjective as a map of sheaves), which implies that $p_{\diag}=p_{\wt{X}}$.

The fact that $\pi_{\diag}\ne\pi_{\wt{X}}$ is equivalent to $W_{\wt{X}}\ne\cE_{K}\times_{Y(\Gamma)_{K}}\cdots\times_{Y(\Gamma)_{K}}\cE_{K}$. There are many ways of seeing this -- we will soon see that the monodromy representations of their relative $H^{1}_{\et}$ are different. A more elementary way of seeing the difference is to observe that the descent data for the two schemes are different for the \'etale covering $\alpha:\wt{Y}\rar Y(\Gamma)_{K}$. Let us fix an isomorphism
\[\iota:\wt{Y}\coprod\cdots\coprod\wt{Y}\cong\wt{Y}\times_{Y(\Gamma)_{K}}\wt{Y},\]such that the two projection maps $p_{1},p_{2}:\wt{Y}\times_{Y(\Gamma)_{K}}\wt{Y}\rar \wt{Y}$ are identified with 
\[p_{1}\circ\iota:\wt{Y}\coprod\cdots\coprod\wt{Y}\xrar{(\id,\cdots,\id)}\wt{Y},\]
\[p_{2}\circ\iota:\wt{Y}\coprod\cdots\coprod\wt{Y}\xrar{(\sigma_{1},\cdots,\sigma_{r})}\wt{Y},\]where $\Gal(\wt{Y}/Y(\Gamma)_{K})=\lbrace\sigma_{1}=\id,\sigma_{2},\cdots,\sigma_{r}\rbrace$. Note that $\Gal(\wt{Y}/Y(\Gamma)_{K})\subset S_{r}$ where one can identify $\sigma\in\Gal(\wt{Y}/Y(\Gamma)_{K})$ with the permutation of the components
\[\wt{Y}\coprod\cdots\coprod\wt{Y}\riso\wt{Y}\times_{Y(\Gamma)_{K}}\wt{Y}\xrar{(\sigma,\id)}\wt{Y}\times_{Y(\Gamma)_{K}}\wt{Y}\liso\wt{Y}\coprod\cdots\coprod\wt{Y}.\]
The descent datum for $\cE_{K}\times_{Y(\Gamma)_{K}}\cdots\times_{Y(\Gamma)_{K}}\cE_{K}$ for the covering $\wt{Y}\rar Y(\Gamma)_{K}$ is given by \[\left(\alpha^{*}\cE_{K}\times_{\wt{Y}}\cdots\times_{\wt{Y}}\alpha^{*}\cE_{K}\right)\coprod\cdots\coprod\left(\alpha^{*}\cE_{K}\times_{\wt{Y}}\cdots\times_{\wt{Y}}\alpha^{*}\cE_{K}\right)\]\[\xrar{(\id,\cdots,\id)}\left(\alpha^{*}\cE_{K}\times_{\wt{Y}}\cdots\times_{\wt{Y}}\alpha^{*}\cE_{K}\right)\coprod\cdots\coprod\left(\alpha^{*}\cE_{K}\times_{\wt{Y}}\cdots\times_{\wt{Y}}\alpha^{*}\cE_{K}\right),\]whereas the descent datum for $W_{\wt{X}}$ for the covering $\wt{Y}\rar Y(\Gamma)_{K}$ is given by\[\left(\alpha^{*}\cE_{K}\times_{\wt{Y}}\cdots\times_{\wt{Y}}\alpha^{*}\cE_{K}\right)\coprod\cdots\coprod\left(\alpha^{*}\cE_{K}\times_{\wt{Y}}\cdots\times_{\wt{Y}}\alpha^{*}\cE_{K}\right)\]\[\xrar{(f,\cdots,f)}\left(\alpha^{*}\cE_{K}\times_{\wt{Y}}\cdots\times_{\wt{Y}}\alpha^{*}\cE_{K}\right)\coprod\cdots\coprod\left(\alpha^{*}\cE_{K}\times_{\wt{Y}}\cdots\times_{\wt{Y}}\alpha^{*}\cE_{K}\right),\]where \[f:\alpha^{*}\cE_{K}\times_{\wt{Y}}\cdots\times_{\wt{Y}}\alpha^{*}\cE_{K}\rar\left(\alpha^{*}\cE_{K}\times_{\wt{Y}}\cdots\times_{\wt{Y}}\alpha^{*}\cE_{K}\right)\coprod\cdots\coprod\left(\alpha^{*}\cE_{K}\times_{\wt{Y}}\cdots\times_{\wt{Y}}\alpha^{*}\cE_{K}\right),\]is $f=(\sigma_{1},\cdots\sigma_{r})$, where $\sigma_{i}$ is the natural map corresponding to the permutation that it represents. As there is no $\wt{Y}$-automorphism that intertwines the two descent data, the two descend to two non-isomorphism $Y(\Gamma)_{K}$-schemes, as desired. 
\end{proof}
\begin{prop2}
The family of abelian varieties $W_{\wt{X}}$ on $Y(\Gamma)_{K}$ has the following properties.
\begin{enumerate}
\item Given a field extension $K'/K$ and a point $x\in Y(\Gamma)(K')$, the fiber $(W_{\wt{X}})_{x}$ is the abelian variety over $K'$ given by
\[(W_{\wt{X}})_{x}=\tR_{\wt{X}_{x}/K'}((E_{x})_{\wt{X}_{x}}),\]where $E_{x}$ is the elliptic curve corresponding to $x$, and $\wt{X}_{x}$ is the \'etale $K'$-algebra of degree $r$ given by the fiber of $\alpha:\wt{X}\rar X(\Gamma)_{K}$ over $x$. 
\item The family $u:W_{\wt{X}}\rar Y(\Gamma)_{K}$ is $\Gal(\wt{Y}/Y(\Gamma)_{K})$-equivariant, where $\Gal(\wt{Y}/Y(\Gamma)_{K})$ acts trivially on $Y(\Gamma)_{K}$.
\end{enumerate}
\end{prop2}

\begin{proof}
(1) follows directly from the fact that the Weil restriction of schemes is compatible with base-change. The action of $\Gal(\wt{Y}/Y(\Gamma)_{K})$ on $\wt{Y}$ and on $\alpha^{*}\cE_{K}$ gives, by functoriality of the Weil restriction, the action of $\Gal(\wt{Y}/Y(\Gamma)_{K})$ on $W_{\wt{X}}$, fixing $Y(\Gamma)_{K}$ on the base, from which (2) follows.
\end{proof}
The twisted Kuga--Sato variety $W_{\wt{X}}$ gives rise to the ``geometric local systems'' on the modular curve, which will be related to the uniformizing logarithmic Higgs bundles.
%
\begin{defn2}\label{LocSysConstruction}
Suppose that $\alpha:\wt{X}\rar X(\Gamma)_{K}$ is an abelian \'etale Galois cover, and that the exponent of the abelian group $\Gal(\wt{X}/X(\Gamma)_{K})$ is $n$. For a character $\chi:\Gal(\wt{X}/X(\Gamma)_{K})\rar \bZ[\zeta_{n}]^{\times}$, we define the following ``geometric local systems'':
\begin{itemize}
\item for a complex embedding $\sigma:K\hrar\bC$, we define a varitaion of polarized pure $\bZ[\zeta_{n}]$-Hodge structures of weight $1$ and rank $2$ on $Y(\Gamma)$ (as a Riemann surface),
\[\rho_{\wt{X},\chi,H}:=\sH^{1}_{B}(W_{\wt{X}}\times_{K,\sigma}\bC/Y(\Gamma),\bZ[\zeta_{n}])[\chi],\]where $\sH^{1}_{B}(W_{\wt{X}}\times_{K,\sigma}\bC/Y(\Gamma),\bZ[\zeta_{n}])$ is the relative first Betti cohomology with $\bZ[\zeta_{n}]$-coefficients, 
and  $[\chi]$ means the $\chi$-isotypic part of the action of $\Gal(\wt{X}/X(\Gamma)_{K})$ on $\sH_{B}^{1}(W_{\wt{X}}\times_{K,\sigma}\bC/Y(\Gamma))$;
\item for $(p,Nn)=1$, we define an \'etale $\bZ_{p}[\zeta_{n}]$-local system of rank $2$ over $Y(\Gamma)_{K}$,
\[\rho_{\wt{X},\chi,p}:=R^{1}u_{\et,*}{\bZ_{p}[\zeta_{n}]}[\chi],\]where $u:W_{\wt{X}}\rar Y(\Gamma)_{K}$ is the natural map, and $[\chi]$ means the $\iota\circ\chi$-isotypic part of the action of $\Gal(\wt{X}/X(\Gamma)_{K})$ on $R^{1}u_{\et,*}{\bZ_{p}[\zeta_{n}]}$, where $\iota:\bZ[\zeta_{n}]\hrar\bZ_{p}[\zeta_{n}]$ is the natural embedding.
\end{itemize}
\end{defn2}
As these local systems come from geometry, the following properties are immediate.
\begin{itemize}
\item The Betti-\'etale comparison isomorphism holds: namely, 
\[\rho_{\wt{X},\chi,p}\rvert_{\pi_{1,\et}(Y(\Gamma)_{K}\times_{K,\sigma}\bC,*)}\cong\wh{\iota\circ\rho_{\wt{X},\chi,H}}:\wh{\pi_{1}(Y(\Gamma),*)}\rar\GL_{2}(\bZ_{p}[\zeta_{n}]).\]Here, the left hand side is the restriction of the $\bZ_{p}[\zeta_{n}]$-\'etale local system $\rho_{\wt{X},\chi,p}$ to the geometric fundamental group $\pi_{1,\et}(Y(\Gamma)_{K}\times_{K,\sigma}\bC,*)$, which is naturally isomorphic to the profinite completion of the topological fundamental group $\pi_{1}(Y(\Gamma),*)$, and the right hand side is the profinite completion of the topological local system $\iota\circ\rho_{\wt{X},\chi,H}:\pi_{1}(Y(\Gamma),*)\rar\GL_{2}(\bZ_{p}[\zeta_{n}])$.
\item For $(p,Nn)=1$, the \'etale local system $\rho_{\wt{X},\chi,p}$ extends to an \'etale local system over an integral model $\fY$ over $\cO_{K,S}$ of $Y(\Gamma)_{K}$ for a finite set of primes $S$ of $K$ including the primes above $p$. Furthermore, at every place $\fp$ of $K$ above $p$, $\rho_{\wt{X},\chi,p}\rvert_{\pi_{1,\et}(Y(\Gamma)_{K_{\fp}},*)}$ is a de Rham local system\footnote{For the definition of de Rham local systems, see \cite{Scholze}.}. 
\end{itemize}

We now specialize the above situation to the case of theta characteristics. Let $K/\bQ$ be a number field, and $\nu$ be a theta characteristic with Kodaira--Spencer data over $X(\Gamma)_{K}$, so that $\nu=\omega\otimes L$ for a $2$-torsion line bundle $L$. Let $\alpha:\wt{X}_{L}\rar X(\Gamma)_{K}$ be the Galois double cover corresponding to the $2$-torsion line bundle $L$ via Lemma \ref{DoubleCoverLem}. We will show that the $-1$-isotypic part of the geometric local systems, namely the variation of polarized pure $\bZ$-Hodge structures $\rho_{\wt{X}_{L},-1,\sigma,H}$ for ${\sigma:K\hrar\bC}$, and the \'etale $\bZ_{p}$-local system $\rho_{\wt{X}_{L},-1,p}$ for ${(p,2N)=1}$, correspond to the uniformizing logarithmic Higgs bundle $(E_{\nu}=\nu\oplus \nu^{-1},\theta_{\nu})$ via the complex and $p$-adic nonabelian Hodge correspondences, respectively. We will use the simplified notation
\[\rho_{\nu,\sigma}:=\rho_{\wt{X}_{L},-1,\sigma,H},\quad \rho_{\nu,p}:=\rho_{\wt{X}_{L},-1,p}.\]
\begin{rmk2}\label{UnderlyingRemark}
It should be noted that the isomorphism class of the uniformizing logarithmic Higgs bundle $(E_{\nu},\theta_{\nu})$ only depends on the underlying theta characteristic $\nu$ and does not depend on the Kodaira--Spencer isomorphism $\KS_{\nu}$. Indeed, $(E_{\nu},\theta_{\nu})\cong (E_{\nu},t\theta_{\nu})$ as Higgs bundles for any nonzero scalar $t$. Thus, one should only expect that $(E_{\nu},\theta_{\nu})$ may only determine the information that only depends on the underlying line bundle of the $2$-torsion line bundle $L$. 
\end{rmk2}
\subsection{Variation of Hodge structures attached to $(E_{\nu},\theta_{\nu})$}
For the variation of Hodge structures $\rho_{\nu,\sigma}$, we aim to show the following.
\begin{thm2}\label{UniquenessHodge} For a complex embedding $\sigma:K\hrar\bC$, $\rho_{\nu,\sigma}$ is the unique, up to the shift of indices, variation of Hodge structures over $Y(\Gamma)$ where the associated graded of its canonical extension is isomorphic to $(E_{\nu},\theta_{\nu})\times_{K,\sigma}\bC$.
\end{thm2}
\begin{rmk2}This is indeed consistent with Remark \ref{UnderlyingRemark}, as $(\wt{X}_{L})_{\bC}$ only depends on the underlying line bundle of $L$; namely, $\bC^{\times}=(\bC^{\times})^{2}$.
\end{rmk2}
\begin{proof}
Let us use the notation $(-)_{\sigma}$ for the shorthand of $(-)\times_{K,\sigma}\bC$. In the proof of Theorem \ref{MainTheta}, we have already seen that $(E_{\nu,\sigma},\theta_{\nu,\sigma})$ is a stable Higgs bundle. Furthermore, it is clear that $(E_{\nu,\sigma},\theta_{\nu,\sigma})\cong(E_{\nu,\sigma},t\theta_{\nu,\sigma})$ for any $t\in\bC^{\times}$. Therefore, by \cite[Lemma 4.1]{SimpsonTensor}, the local system corresponding to $(E_{\nu,\sigma},\theta_{\nu,\sigma})$ comes from a complex variation of Hodge structures, which is unique up to the shift of indices.

On the other hand, the variation of polarized pure $\bZ$-Hodge structures $\rho_{\nu,\sigma}$ has the underlying vector bundle $\sH_{\nu,\sigma}:=\rho_{\nu,\sigma}\otimes\cO_{Y(\Gamma)}$ isomorphic to the relative de Rham cohomology
\[\sH_{\nu,\sigma}\cong\sH^{1}_{\dR}(W_{\wt{X},\sigma}/Y(\Gamma))[\chi_{\wt{X}}],\]where $\chi_{\wt{X}}:\Gal(\wt{X}/X(\Gamma)_{K})\rar\lbrace\pm1\rbrace\hrar\bZ^{\times}$ is the nontrivial character. Note that the local system underlying the variation of Hodge structures $\sH^{1}_{B}(W_{\wt{X},\sigma}/Y(\Gamma),\bZ)$ is, as the representation of $\pi_{1}(Y(\Gamma),*)$, isomorphic to $\Ind_{\pi_{1}(\wt{Y},*)}^{\pi_{1}(Y(\Gamma),*)}\Res_{\pi_{1}(Y(\Gamma),*)}^{\pi_{1}(\wt{Y},*)}\rho_{E}$, where $\rho_{E}$ is the local system underlying the variation of Hodge structures $\sH_{B}^{1}(\cE/Y(\Gamma),\bZ)$. The fact that $\bH$ is the classifying space of pure polarized $\bZ$-Hodge structures of weight $1$ and rank $2$ implies that, as $Y(\Gamma)\cong\bH/\Gamma$, $\rho_{E}$ is isomorphic to the representation $\pi_{1}(Y(\Gamma),*)\cong\Gamma\hrar\SL_{2}(\bZ)\rar\GL_{2}(\bZ)$.

It is a general fact that, if $H\le G$ is a finite index subgroup, given a representation $\rho$ of $G$,
\[\Ind_{H}^{G}\Res_{G}^{H}\rho\cong \rho\otimes S_{G/H},\]where $S_{G/H}:G\rar\GL_{[G:H]}(\bZ)$ is the left regular representation. Therefore, the local system underlying the VHS $\sH_{B}^{1}(W_{\wt{X},\sigma}/Y(\Gamma),\bZ)$ is isomorphic to $\rho\oplus \rho\otimes\chi_{\wt{X}}$. Thus, the local system underlying $\rho_{\nu,\sigma}$ is isomorphic to $\rho\otimes\chi_{\wt{X}}$. By Proposition \ref{DoublecoverC}, it follows that the canonical extension of $\rho_{\nu,\sigma}$ has the associated Hodge bundle equal to $(E_{\nu,\sigma},\theta_{\nu,\sigma})$. This proves the Theorem.
\end{proof}
\subsection{Crystalline local systems attached to $(E_{\nu},\theta_{\nu})$}
Let $(p,2N\disc(K/\bQ))=1$, and let $\rho_{\nu,p}:\pi_{1}(Y(\Gamma)_{K},*)\rar\GL_{2}(\bZ_{p})$ be the $\bZ_{p}$-\'etale local system defined above. Let $\fp$ be a prime of $K$ lying over $p$ (by the assumption on $p$, $\fp$ is unramified over $p$), and define $\rho_{\nu,p,\fp}:\pi_{1}(Y(\Gamma)_{K_{\fp}},*)\rar\GL_{2}(\bZ_{p})$ to be the restriction of $\rho_{\nu,p}$. In this subsection, we aim to describe the relation between $\rho_{\nu,p,\fp}$ and the uniformizing logarithmic Higgs bundle $(E_{\nu},\theta_{\nu})$ in the optic of the $p$-adic nonabelian Hodge correspondence as developed in \cite{LSZ2} (or its generalization to logarithmic Higgs bundles in \cite{LSYZ}). In what follows, we will freely use the terms from \cite{LSZ2}, \cite{LSYZ} and \cite{TT}.

As Remark \ref{UnderlyingRemark} suggests, the $p$-adic nonabelian Hodge correspondence determines only a part of $\rho_{\nu,p,\fp}$. To explain this, we introduce the following definition.
\begin{defn2}Let $F$ be a field extension of $\bQ_{p}$, where $(p,2N)=1$ (so that $X(\Gamma)_{F}$ has a canonical integral model over $\cO_{F}$). We say that a theta characteristic with Kodaira--Spencer data $(\nu,\KS_{\nu})$ \emph{admits an integral model} if there is a line bundle $\nu_{\cO_{F}}$ on $X(\Gamma)_{\cO_{F}}$ and an isomorphism $\KS_{\nu_{\cO_{F}}}:\nu_{\cO_{F}}^{\otimes 2}\riso\Omega_{X(\Gamma)_{\cO_{F}}/\cO_{F}}^{1}(D)$ whose generic fiber coincides with $(\nu,\KS_{\nu})$. Similarly, we say that a $2$-torsion line bundle $L^{\otimes2}\riso\cO_{X(\Gamma)_{F}}$ \emph{admits an integral model} if the flat extension $\cL\in\Jac(X(\Gamma)_{\cO_{F}})[2](\cO_{F})$ (which uniquely exists as $\Jac(X(\Gamma)_{\cO_{F}})[2]$ is a finite flat group scheme over $\cO_{F}$) admits an isomorphism $\cL^{\otimes2}\riso\cO_{X(\Gamma)_{\cO_{F}}}$ whose generic fiber coincides with the isomorphism $L^{\otimes2}\riso\cO_{X(\Gamma)_{F}}$.
\end{defn2}
\begin{lem2}Let $F$ be a discretely valued field extension of $\bQ_{p}$, with $(p,2N)=1$, and let $\pi\in F$ be a uniformizer. Let $(\nu,\KS_{\nu})$ be a theta characteristic with Kodaira--Spencer data on $X(\Gamma)_{F}$. Then, $(\nu,\KS_{\nu})$ admits an integral model if and only if $(\nu,\pi\KS_{\nu})$ does not admit an integral model. Equivalently, if $\iota:L^{\otimes2}\riso\cO_{X(\Gamma)_{F}}$ is a $2$-torsion line bundle, then $(L,\iota)$ admits an integral model if and only if $(L,\pi\iota)$ does not admit and integral model. In particular, the property ``admitting an integral model'' is insensitive to an unramified base change of $F$.
\end{lem2}
\begin{proof}
Suppose that $(\nu,\KS_{\nu})$ admits an integral model. If $(\nu,t\KS_{\nu})$, for $t\in F^{\times}/(F^{\times})^{2}$, admits an integral model, then the integral Kodaira--Spencer isomorphisms are scalar multiples of each other. This implies that $t\in\cO_{F}^{\times}/(\cO_{F}^{\times})^{2}$. Conversely, if $t\in\cO_{F}^{\times}/(\cO_{F}^{\times})^{2}$, then clearly $(\nu,t\KS_{\nu})$ admits an integral model. Note that there is a short exact sequence
\[1\rar(\cO_{F}^{\times})/(\cO_{F}^{\times})^{2}\rar F^{\times}/(F^{\times})^{2}\xrar{v_{F}}\lbrace\pm1\rbrace\rar1,\]where $v_{F}$ is the normalized valuation on $F$ (i.e. $v_{F}(\pi)=1$). Therefore, $(\nu,\KS_{\nu})$ admits an integral model if and only if $(\nu,\pi\KS_{\nu})$ does not admit an integral model.

Let $F'/F$ be an unramified field extension. If $(\nu,\KS_{\nu})$ admits an integral model, then clearly $(\nu_{F'},\KS_{\nu_{F'}})$ admits an integral model via pullback. If $(\nu,\KS_{\nu})$ does not admit an integral model, then $(\nu,\pi\KS_{\nu})$ admits an integral model, so $(\nu_{F'},\pi\KS_{\nu_{F'}})$ admits an integral model, which implies that $(\nu_{F'},\pi\KS_{\nu_{F'}})$ does not admit an integral model. This finishes the proof.
\end{proof}
\begin{lem2}
Let $F$ be a discretely valued field extension of $\bQ_{p}$ with $(p,2N)=1$ such that the residue field of $F$ is algebraically closed (e.g. the maximal unramified extension of $\bQ_{p}$). Then, given a theta characteristic $\nu$, there are exactly two isomorphism classes of theta characteristic with Kodaira--Spencer data whose underlying theta characteristic is $\nu$: one that admits an integral model and one that does not admit an integral model. The same applies to $2$-torsion line bundles.
\end{lem2}
\begin{proof}
This follows from $\cO_{F}^{\times}=(\cO_{F}^{\times})^{2}$, which is a consequence of Hensel's lemma.
\end{proof}
We first prove that the \'etale $\bZ_{p}$-local system associated with a theta characteristic with Kodaira--Spencer data admitting an integral model is \emph{crystalline}.
\begin{prop2}If $(\nu,\KS_{\nu})$ admits an integral model, the \'etale local system $\rho_{\nu,p,\fp}$ is a crystalline $\bZ_{p}$-local system\footnote{For the definition of crystalline $\bZ_{p}$-local systems, see \cite[Definition 3.10]{TT}.}.
\end{prop2}
\begin{proof}By \cite{DR}, the universal generalized elliptic curve $\ov{\cE}_{K}\rar X(\Gamma)_{K}$ has a natural model over $\cO_{K_{\fp}}$, denoted as $\ov{f}:\ov{\cE}_{\cO_{K_{\fp}}}\rar X(\Gamma)_{\cO_{K_{\fp}}}$, which is smooth over $Y(\Gamma)_{\cO_{K_{\fp}}}$. Furthermore, if we denote $k_{\fp}$ by the residue field of $\fp$, then $X(\Gamma)_{k_{\fp}}$ is a smooth curve. 

The same definition of $\omega$,
\[\omega:=f_{*}\Omega^{1}_{\cE_{\cO_{K_{\fp}}}/Y(\Gamma)_{\cO_{K_{\fp}}}},\]gives rise to a theta characteristic on $Y(\Gamma)_{\cO_{K_{\fp}}}$, as the Kodaira--Spencer isomorphism holds on the integral level \cite[A.1.3.17]{Katz}. There is the unique \emph{canonical extension} on the integral level,
\[\omega^{\can}=\ov{f}_{*}\Omega^{1}_{\ov{\cE}_{\cO_{K_{\fp}}}/X(\Gamma)_{\cO_{K_{\fp}}}}(\log \infty_{f}):=\ov{f}_{*}\left(\Omega^{1}_{\ov{\cE}_{\cO _{K_{\fp}}}/\cO_{K_{\fp}}}(\ov{f}^{-1}(D))/\ov{f}^{*}(\Omega^{1}_{X(\Gamma)_{\cO_{K_{\fp}}}/\cO_{K_{\fp}}}(D))\right),\]which satisfies the Kodaira--Spencer isomorphism\footnote{We were unable to locate a literature that states the log version of the Kodaira--Spencer isomorphism on the integral modular curve. A much more general version of the log Kodaira--Spencer isomorphism on the integral level is proved in \cite[Proposition 6.9]{Lan}, which contains the statements that we would like for the modular curves.}
\[(\omega^{\can})^{\otimes2}\riso\Omega^{1}_{X(\Gamma)_{\cO_{K_{\fp}}}/\cO_{K_{\fp}}}(D).\]Moreover, the canonical extension $\omega^{\can}$ arises as the first Hodge filtration of the log-de Rham cohomology bundle $R^{1}\ov{f}_{\log\dR,*}(\ov{\cE}_{\cO_{K_{\fp}}}/X(\Gamma)_{\cO_{K_{\fp}}})$, where the log-structures for $\ov{\cE}_{\cO_{K_{\fp}}}$ and $X(\Gamma)_{\cO_{K_{\fp}}}$ are given by $\ov{f}^{-1}(D)$ and $D$, respectively. We will omit the superscript $^{\can}$ from now on.


Let $\cL\in\Jac(X(\Gamma)_{\cO_{K_{\fp}}})[2](\cO_{K_{\fp}})$ with $\cL^{\otimes2}\riso\cO_{X(\Gamma)_{K_{\fp}}}$ be an integral model of the $2$-torsion line bundle $L$. By applying the same construction, we obtain the double Galois cover $\alpha:\wt{X}_{\cO_{K_{\fp}}}\rar X(\Gamma)_{\cO_{K_{\fp}}}$ of degree $2$. Accordingly, we obtain $u:\ov{W}_{\wt{X}_{\cO_{K_{\fp}}}}\rar X(\Gamma)_{\cO_{K_{\fp}}}$.
Using this integral model, one can also extend a logarithmic Higgs sheaf $(E_{\nu},\theta_{\nu})=(E_{\omega},\theta_{\omega})\otimes(L,0)$ over $X(\Gamma)_{K_{\fp}}$ to $(E_{\omega},\theta_{\omega})\otimes(\cL,0)$ over $X(\Gamma)_{\cO_{K_{\fp}}}$.

Since $u:W_{\wt{X}_{\cO_{K_{\fp}}}}\rar Y(\Gamma)_{\cO_{K_{\fp}}}$ is smooth and proper, \cite[Proposition 5.4]{TT} implies that the relative crystalline cohomology $\sE:=R^{1}u_{\cris,*}\cO$ gives rise to a convergent $F$-isocrystal on $Y(\Gamma)_{k_{\fp}}$ (which is in fact overconvergent). From the relative crystalline comparison theorem, \cite[Theorem 5.5]{TT}, it follows that $\rho_{\nu,p,\fp}$ is a crystalline $\bZ_{p}$-local system, and is associated to $\sE$.
\end{proof}
\begin{rmk2}
This also implies that, for $(\nu,\KS_{\nu})$ not admitting an integral model, $\rho_{\nu,p,\fp}$ is ``potentially crystalline'' in some sense, but we are not aware of a good formalism for such local systems.
\end{rmk2}
Now we turn to the $p$-adic nonabelian Hodge correspondence as developed by \cite{LSZ2} and \cite{LSYZ}, which associates crystalline local systems to \emph{periodic (logarithmic) Higgs-de Rham flows}. For their definition, see \cite[Definition 5.2]{LSZ2}\footnote{Strictly speaking, \cite{LSZ2} only concerns the non-log case, but the same definition works for the log case using the inverse Cartier transform in the log case as considered in \cite[\S5, \S6]{LSYZ}.}. Combined with the $p$-adic Riemann--Hilbert correspondence of \cite{Faltings}\footnote{A crystalline $\bZ_{p}$-local system in the sense of \cite{Faltings} is a crystalline $\bZ_{p}$-local system in the sense of \cite{TT} by \cite[Proposition 3.21]{TT}.}, we state the version of the $p$-adic nonabelian Hodge correspondence we will use.
\begin{thm2}[($p$-adic nonabelian Hodge correspondence)]\label{HdR2CrysLocSys}Let $k=\ov{\bF}_{p}$ and $F=W(k)[1/p]$. Let $X_{\infty}$ be a smooth projective curve over $W(k)$ and $D_{\infty}\subset X_{\infty}$ be a relative effective Cartier divisor over $W(k)$. Then, there is a fully faithful covariant functor from the category of $1$-periodic logarithmic Higgs-de Rham flow over $X_{\infty}$ to the category of crystalline $\bZ_{p}$-local systems over $Y_{\infty,F}:=(X_{\infty}-D_{\infty})_{F}$ (with respect to the integral model $Y_{\infty}=X_{\infty}-D_{\infty}$) with the Hodge--Tate weights in $[0,p-2]$. This functor preserves the rank.
\end{thm2}
Using this, we will prove the $p$-adic analogue of Theorem \ref{UniquenessHodge}.
\begin{thm2}\label{UniquenessP}
Let $\nu$ be a theta characteristic with a Kodaira--Spencer data over $K$ that admits an integral model. Then, $\rho_{\nu,p,\fp}\rvert_{\pi_{1,\et}(Y(\Gamma)_{K_{\fp}^{\nr}})}$ depends only on the underlying theta characteristic, where $K_{\fp}^{\nr}$ is the maximal unramified extension of $K_{\fp}$. It is the crystalline $\bZ_{p}$-local system associated to a unique filtered convergent $F$-isocrystal on $Y(\Gamma)_{\ov{k}_{\fp}}$ whose associated graded is isomorphic to $(E_{\nu},\theta_{\nu})\times_{K}K_{\fp}^{\nr}$.
\end{thm2}
\begin{proof}
For simplicity, we let $X_{\infty}=X(\Gamma)_{W(\ov{k}_{\fp})}$, with the log structures coming from $D$.
Thanks to the $p$-adic nonabelian Hodge correspondence (Theorem \ref{HdR2CrysLocSys}), the Theorem will follow if we show that there exists, up to a shift in filtration, a unique $1$-periodic logarithmic Higgs-de Rham flow over $X_{\infty}$ whose underlying Higgs bundle is $(E_{\nu,W(\ov{k}_{\fp})},\theta_{\nu,W(\ov{k}_{\fp})})=(E_{\omega,W(\ov{k}_{\fp})},\theta_{\omega,W(\ov{k}_{\fp})})\otimes(\cL_{W(\ov{k}_{\fp})},0)$. As the rank of the Higgs bundle is $2$, the Higgs field is nonzero, and $\nu\ne\nu^{-1}$, the filtration is unique up to shift of indices. Therefore, the vector bundle with flat connection part of the $1$-periodic Higgs-de Rham flow is determined by the inverse Cartier transform of \cite[\S4][LSZ2] and \cite[\S5, \S6]{LSYZ}. This finishes the proof.
\end{proof}
\begin{rmk2}
As \cite[Theorem 1.2]{LSYZ} suggests, the results of this section suggests that $X(\Gamma)_{W(\ov{k}_{\fp})}$ is in some sense a canonical lifting of $X(\Gamma)_{\ov{k}_{\fp}}$ associated to the uniformizing logarithmic Higgs bundle $(E_{\nu, k_{\fp}},\theta_{\nu,k_{\fp}})$.
\end{rmk2}
\section{Twisted period map to a Siegel modular threefold}
From now on, we assume for simplicity\footnote{This is merely for the simplicity of the moduli problem that the corresponding level structure represents.} that $\Gamma$ is either $\Gamma_{1}(N)$ or $\Gamma(N)$. Recall that $W_{\wt{X}}\rar Y(\Gamma)_{K}$ is a family of abelian surfaces that is different from the square of the universal elliptic curve, but they become isomorphic after an \'etale base-change to $\wt{Y}$:
\[W_{\wt{X}}\not\cong \cE_{K}\times_{Y(\Gamma)_{K}}\cE_{K},\quad\alpha^{*}W_{\wt{X}}\cong\alpha^{*}\cE_{K}\times_{\wt{Y}}\alpha^{*}\cE_{K}.\]In this case, it turns out that, from the universal level structure on $\cE/Y(\Gamma)$, one can construct a certain natural level structure on $W_{\wt{X}}/Y(\Gamma)_{K}$, which is a twisted version of the natural level structure on $\cE_{K}\times_{Y(\Gamma)_{K}}\cE_{K}$. This implies that $Y(\Gamma)_{K}$ admits a twisted period map into the moduli space of abelian surfaces with a level structure that is different from the usual ``diagonal embedding.'' Under this diagonal embedding, we see that a Siegel modular form restricts to a noncongruence modular from of level associated with the theta characteristic $\nu$.

To construct the twisted level structure, we first describe $W_{\wt{X}}$ as a variety over $\wt{Y}$ with a descent datum.
\begin{prop2}\label{DescentDatum}Let $\sigma:\wt{Y}\rar\wt{Y}$ be the nontrivial element of $\Gal(\wt{Y}/Y(\Gamma)_{K})$. Let $\lambda:\wt{Y}\times_{Y(\Gamma)_{K}}\wt{Y}\riso\wt{Y}\coprod\wt{Y}$ be an isomorphism of $\wt{Y}$-schemes such that the following diagram commutes.
\[\xymatrix{\wt{Y}\times_{Y(\Gamma)_{K}}\wt{Y}\ar[rd]_-{\pr_{2}}\ar[rr]^-{\lambda}&&\wt{Y}\coprod\wt{Y}\ar[ld]^-{\id\coprod\sigma}\\ &\wt{Y} &}\]Let $\lambda:W_{\wt{X}}\times_{Y(\Gamma)_{K}}\wt{Y}\riso\alpha^{*}\cE_{K}\times_{\wt{Y}}\alpha^{*}\cE_{K}$ be the natural isomorphism obtained from $\lambda:\wt{Y}\times_{Y(\Gamma)_{K}}\wt{Y}\riso\wt{Y}\coprod\wt{Y}$. Then, the following diagram commutes.
\[\xymatrix{W_{\wt{X}}\times_{Y(\Gamma)_{K}}\wt{Y}\ar[rrr]^-{(\id,\sigma)}\ar[d]_-{\lambda}&&&W_{\wt{X}}\times_{Y(\Gamma)_{K}}\wt{Y}\ar[d]^-{\lambda}\\\alpha^{*}\cE_{K}\times_{\wt{Y}}\alpha^{*}\cE_{K}\ar[rrr]_-{(x,y)\mapsto(\sigma(y),\sigma(x))}&&&\alpha^{*}\cE_{K}\times_{\wt{Y}}\alpha^{*}\cE_{K}}\]
\end{prop2}
\begin{proof}Note that \[\wt{Y}\times_{Y(\Gamma)_{K}}\wt{Y}\xrar{(\id,\sigma)}\wt{Y}\times_{Y(\Gamma)_{K}}\wt{Y},\]after conjugating by $\lambda$, is identified with \[\wt{Y}\coprod\wt{Y}\xrar{x\coprod y\mapsto y\coprod x}\wt{Y}\coprod\wt{Y}.\]By \cite[Proposition A.5.2]{CGP}, the Weil restriction of schemes has a natural isomorphism
\[\tR_{S'/S}(X')\times_{S}T\cong\tR_{T'/T}(X'\times_{S'}T'),\]where $S'$ is a finite locally free $S$-scheme, $X'$ is an $S'$-scheme and $T'=S'\times_{S}T$. Therefore, the isomorphism \[W_{\wt{X}}\times_{Y(\Gamma)_{K}}\wt{Y}=\left(\tR_{\wt{Y}/Y(\Gamma)_{K}}(\alpha^{*}\cE_{K})\right)\times_{Y(\Gamma)_{K}}\wt{Y}\cong \tR_{\wt{Y}\times_{Y(\Gamma)_{K}}\wt{Y}/\wt{Y}}(\alpha^{*}\cE_{K}\times_{\wt{Y}}(\wt{Y}\times_{Y(\Gamma)_{K}}\wt{Y})),\]is natural, where in the rightmost expression, the morphism $\wt{Y}\times_{Y(\Gamma)_{K}}\wt{Y}\rar\wt{Y}$ used in the subscript is the second projection, while the morphism $\wt{Y}\times_{Y(\Gamma)_{K}}\wt{Y}\rar\wt{Y}$ used in the expression in the parenthesis is the first projection. Thus, after conjugating by $\lambda:\wt{Y}\times_{Y(\Gamma)_{K}}\wt{Y}\riso\wt{Y}\coprod\wt{Y}$, this is identified with 
\[\tR_{\wt{Y}\coprod\wt{Y}/\wt{Y}}(\alpha^{*}\cE_{K}\times_{\wt{Y}}(\wt{Y}\coprod\wt{Y})),\]where $\id\coprod\sigma:\wt{Y}\coprod\wt{Y}\rar\wt{Y}$ ($\id\coprod\id:\wt{Y}\coprod\wt{Y}\rar\wt{Y}$, respectively) is used in the subscript (the expression in the parenthesis, repsectively). Therefore, under this identification, the morphism
\[(\id,\sigma):W_{\wt{X}}\times_{Y(\Gamma)_{K}}\wt{Y}\rar W_{\wt{X}}\times_{Y(\Gamma)_{K}}\wt{Y},\]is identified with the morphism
\[\tR_{\wt{Y}\coprod\wt{Y}/\wt{Y}}(\alpha^{*}\cE_{K}\times_{\wt{Y}}(\wt{Y}\coprod\wt{Y}))\rar\tR_{\wt{Y}\coprod\wt{Y}/\wt{Y}}(\alpha^{*}\cE_{K}\times_{\wt{Y}}(\wt{Y}\coprod\wt{Y})),\]where the subscripts are related by the diagram
\[\xymatrix{\wt{Y}\coprod\wt{Y}\ar[rrr]^-{x\coprod y\mapsto y\coprod x}\ar[d]_{\id\coprod\sigma}&&&\wt{Y}\coprod\wt{Y}\ar[d]^-{\id\coprod\sigma}\\\wt{Y}\ar[rrr]_-{\sigma}&&&\wt{Y}}\]
and the expressions in the parentheses are related by the diagram
\[\xymatrix{\wt{Y}\coprod\wt{Y}\ar[rr]^-{x\coprod y\mapsto y\coprod x}\ar[rd]_-{\id\coprod\id}&&\wt{Y}\coprod\wt{Y}\ar[ld]^-{\id\coprod\id}\\&\wt{Y}&}\]
From this, the statement easily follows.
\end{proof} 
We consider a $\Gamma$-level structure on an elliptic scheme $E/S$. In the case of $\Gamma=\Gamma(N)$, it is a pair of sections $P_{1},P_{2}:S\rar E$ that fiberwise generates $E[N]$, and in the case of $\Gamma=\Gamma_{1}(N)$, it is a section $P:S\rar E[N]$ that has exact order $N$. We take the $\Gamma$-level structure on the universal elliptic curve $\cE_{K}/Y(\Gamma)_{K}$ as either $\cP,\cQ:Y(\Gamma)_{K}\rar\cE_{K}[N]$ (in the case of $\Gamma(N)$) or $\cP:Y(\Gamma)_{K}\rar\cE_{K}[N]$ (in the case of $\Gamma_{1}(N)$). 
Using the level structure on $\cE_{K}$, we may define a twisted level structure on $W_{\wt{X}}$ as follows.
\begin{defn2}[($\Gamma(N)^{+}$- and $\Gamma_{1}(N)^{+}$-structures on an abelian surface)]For a principally polarized abelian surface $(A/S,\lambda)$, a  $\Gamma(N)^{+}$-structure is a collection of \'etale-local sections $P_{1},P_{2},P_{3},P_{4}$ of $A[N]$ such that they generate $A[N]$ fiberwise, and two such collections \[P_{1},P_{2},P_{3},P_{4},\quad P_{1}',P_{2}',P_{3}',P_{4}',\]are equivalent if $\lbrace P_{1},P_{2}\rbrace=\lbrace P_{1}',P_{2}'\rbrace$ and $\lbrace P_{3},P_{4}\rbrace=\lbrace P_{3}',P_{4}'\rbrace$ (as unordered sets). 

A  $\Gamma_{1}(N)^{+}$-structure is a collection of \'etale-local sections $P_{1},P_{2}$ of $A[N]$ such that they generate a totally isotropic subspace of $A[N]$, with respect to the Weil pairing induced by $\lambda$, and two such collections \[P_{1},P_{2},\quad P_{1}',P_{2}',\]are equivalent if $\lbrace P_{1},P_{2}\rbrace=\lbrace P_{1}',P_{2}'\rbrace$ (as unordered sets).
\end{defn2}
\begin{rmk2}The moduli space of principally polarized abelian surfaces with $\Gamma(N)^{+}$- or $\Gamma_{1}(N)^{+}$-structures are identified with the arithmetic quotient of the Siegel upper half space by a subgroup of $\Sp_{4}(\bZ)$. More precisely, if the symplectic form corresponds to the matrix
$\begin{psmat}0&1&0&0\\-1&0&0&0\\0&0&0&1\\0&0&-1&0\end{psmat}$, then \[\Gamma(N)^{+}=\Gamma(N)\cdot\left\langle\begin{pmatrix}0&1&0&0\\-1&0&0&0\\0&0&1&0\\0&0&0&1\end{pmatrix}\right\rangle,\qquad\Gamma_{1}(N)^{+}=\Gamma_{1}(N)\cdot\left\langle\begin{pmatrix}0&1&0&0\\-1&0&0&0\\0&0&1&0\\0&0&0&1\end{pmatrix}\right\rangle,\]where $\Gamma(N)$ and $\Gamma_{1}(N)$ are the standard congruence subgroups of $\Sp_{4}(\bZ)$,
\[\Gamma(N)=\lbrace M\in\Sp_{4}(\bZ)~|~M\equiv I_{4}~(\MOD N)\rbrace,\]
\[\Gamma_{1}(N)=\lbrace M\in\Sp_{4}(\bZ)~|~M~(\MOD N)\text{ is upper triangular unipotent}\rbrace.\]
\end{rmk2}
\begin{defn2}[(Twisted level structure on $W_{\wt{X}}$)]
Let $\Gamma^{+}$ be $\Gamma(N)^{+}$ ($\Gamma_{1}(N)^{+}$, respectively) if $\Gamma=\Gamma(N)$ ($\Gamma=\Gamma_{1}(N)$, respectively). 
We define the twisted level structure, a $\Gamma^{+}$-structure on $\alpha^{*}W_{\wt{X}}\cong \alpha^{*}\cE_{K}\times_{\wt{Y}}\alpha^{*}\cE_{K}$, as follows.
\begin{itemize}
\item If $\Gamma=\Gamma(N)$, we consider the sections $\wt{\cP}_{1},\wt{\cP}_{2},\wt{\cQ}_{1},\wt{\cQ}_{2}:\wt{Y}\rar(\alpha^{*}\cE_{K}\times_{\wt{Y}}\alpha^{*}\cE_{K})[N]$, where
\[\wt{\cP}_{1}:=(\alpha^{*}\cP,\alpha^{*}e),\quad\wt{\cP}_{2}:=(\alpha^{*}e,\alpha^{*}\cP),\quad\wt{\cQ}_{1}:=(\alpha^{*}\cQ,\alpha^{*}e),\quad\wt{\cQ}_{2}:=(\alpha^{*}e,\alpha^{*}\cQ).\]The above \'etale-local sections define a $\Gamma(N)^{+}$-structure on $\alpha^{*}W_{\wt{X}}$.
\item If $\Gamma=\Gamma_{1}(N)$, we consider the sections $\wt{\cP}_{1},\wt{\cP}_{2}:\wt{Y}\rar(\alpha^{*}\cE_{K}\times_{\wt{Y}}\alpha^{*}\cE_{K})[N]$, where \[\wt{\cP}_{1}:=(\alpha^{*}\cP,\alpha^{*}e),\quad \wt{\cP}_{2}:=(\alpha^{*}e,\alpha^{*}\cP).\]The above \'etale-local sections define a $\Gamma_{1}(N)^{+}$-structure on $\alpha^{*}W_{\wt{X}}$.
\end{itemize}
\end{defn2}
\begin{lem2}
The twisted level structure on $\alpha^{*}W_{\wt{X}}$, as a $\Gamma^{+}$-level structure, descends into a twisted level structure, again as a $\Gamma^{+}$-level structure, on $W_{\wt{X}}$. Namely, the twisted level structure on $\alpha^{*}W_{\wt{X}}$ is invariant under the automorphism induced by $\sigma:\wt{Y}\rar\wt{Y}$.
\end{lem2}
\begin{proof}
We already know what descent datum $\alpha^{*}W_{\wt{X}}\cong\alpha^{*}\cE_{K}\times_{\wt{Y}}\alpha^{*}\cE_{K}$ has, thanks to Proposition \ref{DescentDatum}. We only need to check that that the $\Gamma^{+}$-level structure is compatible with the descent datum, which is clear as the level structure is indifferent to the switch between $\wt{\cP}_{1}$ and $\wt{\cP}_{2}$ (and also the switch between $\wt{\cQ}_{1}$ and $\wt{\cQ}_{2}$, if $\Gamma=\Gamma(N)$).
\end{proof}
\begin{rmk2}
The twisted period map $\pi_{\wt{X}}:Y(\Gamma)_{K}\rar \cA_{2,\Gamma^{+}}$ is different from the usual ``diagonal'' period map $\pi_{\diag}:Y(\Gamma)_{K}\rar \cA_{2,\Gamma^{+}}$, given by the diagonal morphism between the moduli functors, $E\mapsto E^{2}$. This is simply because the pullbacks of the universal abelian surface over $\cA_{2,\Gamma^{+}}$ by the two period maps are different.

It may first look strange to have a classifying map into a congruence quotient of a Shimura variety even though the starting object is ``noncongruence''. This phenomenon happens because the double cover of congruence quotients $\cA_{2,\Gamma}\rar\cA_{2,\Gamma^{+}}$ of a larger group somehow ``absorbs'' the double cover $\wt{Y}\rar Y(\Gamma)_{K}$. To be more precise, for the diagonal period map, there is a map $Y(\Gamma)_{K}\rar \cA_{2,\Gamma}$ that fills in the diagram
\[\xymatrix{Y(\Gamma)_{K}\coprod Y(\Gamma)_{K}\ar[r]\ar[d]&\cA_{2,\Gamma}\ar[d]\\ Y(\Gamma)_{K}\ar[ru]\ar[r]&\cA_{2,\Gamma^{+}}}\]On the other hand, for the twisted period map $\pi_{\wt{X}}$, the diagonal arrow cannot be filled:
\[\xymatrix{\wt{Y}\ar[r]\ar[d]&\cA_{2,\Gamma}\ar[d]\\ Y(\Gamma)_{K}\ar[r]\ar@{-->}[ru]|-{\times\times\times}&\cA_{2,\Gamma^{+}}}\]
It is interesting to note that we had to use the stacky double cover $\cA_{2,\Gamma}\rar\cA_{2,\Gamma^{+}}$, which seems necessary.
\end{rmk2}
The following suggests that a Siegel modular form restricted along the twisted period map gives rise to a noncongruence modular form.
\begin{prop2}
Let $\omega_{\cA}$ be the Hodge bundle on $\cA_{2,\Gamma^{+}}$, i.e. the automorphic vector bundle corresponding to the standard representation of $\Sp_{4}$, or equivalently, $\omega_{\cA}:=p_{*}\Omega_{\cX/\cA_{2,\Gamma^{+}}}^{1}$, where $p:\cX\rar\cA_{2,\Gamma^{+}}$ is the universal abelian surface. Then, $\pi_{\wt{X}}^{*}\omega_{\cA}=\omega\oplus\nu$, where $\nu$ is the underlying theta characteristic of the theta characteristic with Kodaira--Spencer data associated with the cover $\wt{Y}\rar Y(\Gamma)_{K}$.
\end{prop2}
\begin{proof}
This is a simple consequence of the cohomology and base change theorem, e.g. \cite[Theorem III.12.11]{Hartshorne}.
\end{proof}
This implies that the Siegel modular form of weight $(k,j)$ and level $\Gamma^{+}$, which corresponds to being a section of the automorphic vector bundle $\det^{k}\omega_{\cA}\otimes\Sym^{j}\omega_{\cA}$, restricts along the twisted period map to a section of the vector bundle
\[\omega^{k+j}\otimes\nu^{k}\oplus\omega^{k+j-1}\otimes\nu^{k+1}\oplus\cdots\oplus\omega^{k+1}\otimes\nu^{k+j-1}\oplus\omega^{k}\otimes\nu^{k+j}.\]For example, the restriction of the Siegel modular form of weight $(k,1)$ is a pair of a modular form of weight $2k+1$ and level $\Gamma$ and a \emph{noncongruence} modular form of weight $2k+1$ and level $\Gamma_{\nu}$.
\bibliographystyle{alpha}
\bibliography{Weight1Theta}
\end{document}